\pdfoutput=1
\documentclass[11pt]{article}
\usepackage[margin=1.0in]{geometry}

\usepackage{amsmath,amsfonts,amsthm,amssymb} % Mathematical typesetting and symbols
\usepackage{graphicx} % Figure inclusion
\usepackage{hyperref} % Hyperlinks for citations, references, and URLs
\usepackage[numbers,sort&compress]{natbib} % Citation styling
\usepackage{placeins} % Figure placement control
\usepackage{subcaption} % Captions for subfigures
\usepackage[noend]{algpseudocode} % Algorithm pseudocode
\usepackage{enumerate} % Custom enumerate style

\algrenewtext{Procedure}{} % Using procedure blocks for spacing
\algrenewcommand{\algorithmiccomment}[1]{\hfill // #1} % C-like comments

\newcommand{\set}[2]{\left\{#1\,\middle|\,#2\right\}}

\title{A Public Transit Network Optimization Model for Equitable Access to Social Services}
\author{Adam Rumpf\footnote{Department of Applied Mathematics, Florida Polytechnic University, Lakeland, FL 33805. E-mail: \href{mailto:arumpf@floridapoly.edu}{arumpf@floridapoly.edu}} \, and Hemanshu Kaul\footnote{Department of Applied Mathematics, Illinois Institute of Technology, Chicago, IL 60616. E-mail: \href{mailto:kaul@iit.edu}{kaul@iit.edu}}}
\date{November 11, 2021}

\begin{document}

\maketitle

\begin{abstract}
We present a flexible public transit network design model which optimizes a social access objective while guaranteeing that the system's costs and transit times remain within a preset margin of their current levels. The purpose of the model is to find a set of minor, immediate modifications to an existing bus network that can give more communities access to the chosen services while having a minimal impact on the current network's operator costs and user costs. Design decisions consist of reallocation of existing resources in order to adjust line frequencies and capacities. We present a hybrid tabu search/simulated annealing algorithm for the solution of this optimization-based model. As a case study we apply the model to the problem of improving equity of access to primary health care facilities in the Chicago metropolitan area. The results of the model suggest that it is possible to achieve better primary care access equity through reassignment of existing buses and implementation of express runs, while leaving overall service levels relatively unaffected.
\end{abstract}

%==============================================================================

% Introduction
\section{Introduction}
\label{sec:intro}

The efficient design and operation of a public transit network is one of the most important and challenging tasks facing any urban planner. As such, the \textit{urban transit network design problem} (UTNDP) has been the subject of an immense amount of attention in transportation engineering research for many decades. These problems are generally modeled using an optimization-based framework. Although there is a large amount of variety in the types of design decisions, constraints, and objectives used to define the model, the most common objectives consist of optimizing some combination of operator costs (including fleet size, vehicle maintenance, travel distance, number of stops, profit, etc.) and user costs (including travel time, waiting time, walking distance, number of transfers, vehicle crowding, etc.).

To give a few recent examples of UTNDP studies, Bielli et al.~2002 \cite{bielli2002} proposed a genetic algorithm for bus network design with a wide variety of performance indicators, including minimizing number of required vehicles, number of lines, number of user transfers, user waiting time, and high crowding indices, while maximizing number of served users and vehicle efficiency. Chakroborty and Dwivedi et al.~2002 \cite{chakroborty2002} considered a route design problem with the objective of minimizing a weighted combination of line scores based on minimizing user travel times, maximizing fraction of users able to reach their destination in zero, one, or two transfers, and minimizing fraction of users unable to reach their destination in two or fewer transfers. Fan and Machemehl 2006 \cite{fan2006} used a genetic algorithm to minimize a weighted sum of user costs (consisting of total travel time on both direct and transfer paths), operator costs (consisting of vehicle operation costs), and total unsatisfied transit demand. Bagloee and Ceder 2011 \cite{bagloee2011} studied transit planning decisions with the objective of minimizing a passenger discomfort index consisting of a weighted sum of in-vehicle time, walking time, waiting time, and number of transfers. Cipriani et al.~2012 \cite{cipriani2012} considered a bus route design problem attempting to minimize a weighted combination of user costs (in-vehicle travel time, access time, waiting time, and transfer penalty) and operator costs (bus travel distance and travel time), and a penalty for unsatisfied transit demand. Cheung and Shalaby 2016 \cite{cheung2016} considered rerouting decisions with the goal of minimizing network congestion. Buba and Lee 2018 \cite{buba2018} used a differential evolution algorithm to solve a route design and frequency setting problem with the objective of minimizing a weighted sum of total user travel time and unmet demand. Chu 2018 \cite{chu2018} considered a bus route design and timetabling problem with the objective of minimizing a weighted sum of vehicle travel time, user travel time, and unsatisfied demand.

A relatively small set of operator cost and user cost objectives related to vehicle operation cost, user travel time, and travel demand satisfaction is repeated across most UTNDP studies. This is noted in literature reviews on the topic, such as Guihaire and Hao 2008 \cite{guihaire2008} and Farahani et al.~2013 \cite{farahani2013}, the latter of which specifically points out that other important objectives like environmental costs have gone largely ignored by past studies. It is understandable that these standard operator and user cost objectives are commonly used, since they apply directly to the most immediate and important measures of system performance. Even in a study which includes alternate optimization criteria these objectives are still extremely important to take into consideration. However, such a narrow set of design goals can lead to different communities experiencing extremely disparate levels of access to important public services.

Recently there has been a large amount of interest in studying such questions of differential access to services over the boundaries of a city. Radke and Mu 2000 \cite{radke2000} developed a method for studying household access to social services. Pearce et al.~2006 \cite{pearce2006} studied community access levels to a variety of resources thought to impact community health, including recreation, shopping, education, and health care. Farber et al.~2014 \cite{farber2014} and McDermot et al.~2017 \cite{mcdermot2017} both studied the phenomenon of \textit{food deserts} (areas which lack access to healthy food options) and both found that differential access to food options is largely attributable to differential public transit access as opposed to the geographic locations of facilities. Finally, a large number of recent studies have attempted to measure differential access to primary health care facilities.

Most existing research in this area focuses primarily on developing a health care accessibility metric and using that metric to identify areas with particularly poor access. For example, Delamater 2013 \cite{delamater2013} compared various catchment area-based accessibility metrics for measuring inpatient facility access in the state of Michigan, while McGrail and Humphreys 2014 \cite{mcgrail2014} applied a catchment area-based accessibility metric to measure spatial accessibility of general practitioners throughout Australia. The topic of planning specific actions in order to alleviate these health care access issues has been less well studied. As an example, Gu et al.~2010 \cite{gu2010} considered an optimization-based facility location problem with the goal of placing new primary care facilities to maximize a global health care accessibility metric.

The facility location modeling paradigm is a useful tool for maximizing the public good provided by newly-opened facilities, but creating new facilities is not always an option. The placement of privately-owned facilities is governed by their owners' economic considerations and not necessarily by the needs of the community, and even public facilities are major investments that can generally only be placed in very limited number and in a very limited set of locations. For these reasons we instead turn our attention to an alternative to the facility location modeling paradigm for addressing the accessibility problem.

%%%
\paragraph{Transit Network Design Approach to Accessibility Maximization:}

The primary focus of this study is to consider a problem which is in some ways the inverse of the facility location problem. Rather than changing facility locations in a fixed transit network, we consider changing a transit network containing fixed facility locations. This can improve accessibility levels in the same way that placing facilities can, since reducing travel times makes effective catchment areas larger and brings fast access to more communities. Public transit planning can also be implemented much more immediately since tactical decisions like bus frequency setting can be enacted without the need for permanent changes to the infrastructure.

The main goal of this paper is to present an optimization model-based approach to the public transit design problem that attempts to achieve a social access objective while adhering to the following set of general modeling principles:
\begin{enumerate}
	\item The model should be flexible, allowing the planner to choose any desired social access objective and any assumptions regarding user behavior.
	\item The model must produce solutions that remain at or near the system's current cost and performance level while attempting to optimize the social access objective.
	\item The design decisions should consist of measures which are low-cost, immediate, and easily-implemented.
	\item We will assume that travel related to our social access goal makes up a relatively small proportion of the day-to-day public transit travel volume. In particular we will assume that the capacity of the facilities in question is a much more significant limiting factor than the public transit service capacity.
\end{enumerate}

Principles 2 and 3, in particular, are important for the model to be seriously considered by an urban planner. The main goal of a transit authority is to minimize the operator and user costs resulting from the large volume of day-to-day transit network usage. It would be unreasonable to implement sweeping system changes in service of a relatively small proportion of the overall travel volume, but low-cost, minor, tactical changes that have a minimal effect on most day-to-day traffic might have a better chance of being implemented. As a case study for testing the model formulation presented in this paper, we consider the problem of improving equity of primary health care access in Chicago, Illinois.

%%%
\paragraph{Chicago Metropolitan Area Case Study:}

The study of regional primary access for Chicago residents has been a subject of investigation for the Chicago Department of Public Health. Between 2005 and 2008 a sequence of reports were released to detail their findings \cite{healthcarepuzzle,healthcaresafetynet,healthcareunderserved}. Specifically these reports focused on the so-called \textit{safety net} providers, which consist of a combination of community health centers and free clinics that focus on serving low-income patients.

A number of factors related to patient access to these facilities were studied, including socioeconomic status, ethnicity, gender, age, and geographic location. All reports identified geographic location as one of the biggest contributing factors to access. A large part of this geospatial disparity is due to the locations of the facilities relative to the patient populations. Specifically, Chicago's West Side contains a disproportionately large number of safety net facilities relative to its population, while the North and South Sides contain disproportionately small numbers. In both of these regions estimated population need exceeds estimated facility capacity (as determined by number of available physicians).

\begin{figure}[h]
	\centering
	\includegraphics[height=0.25\textheight]{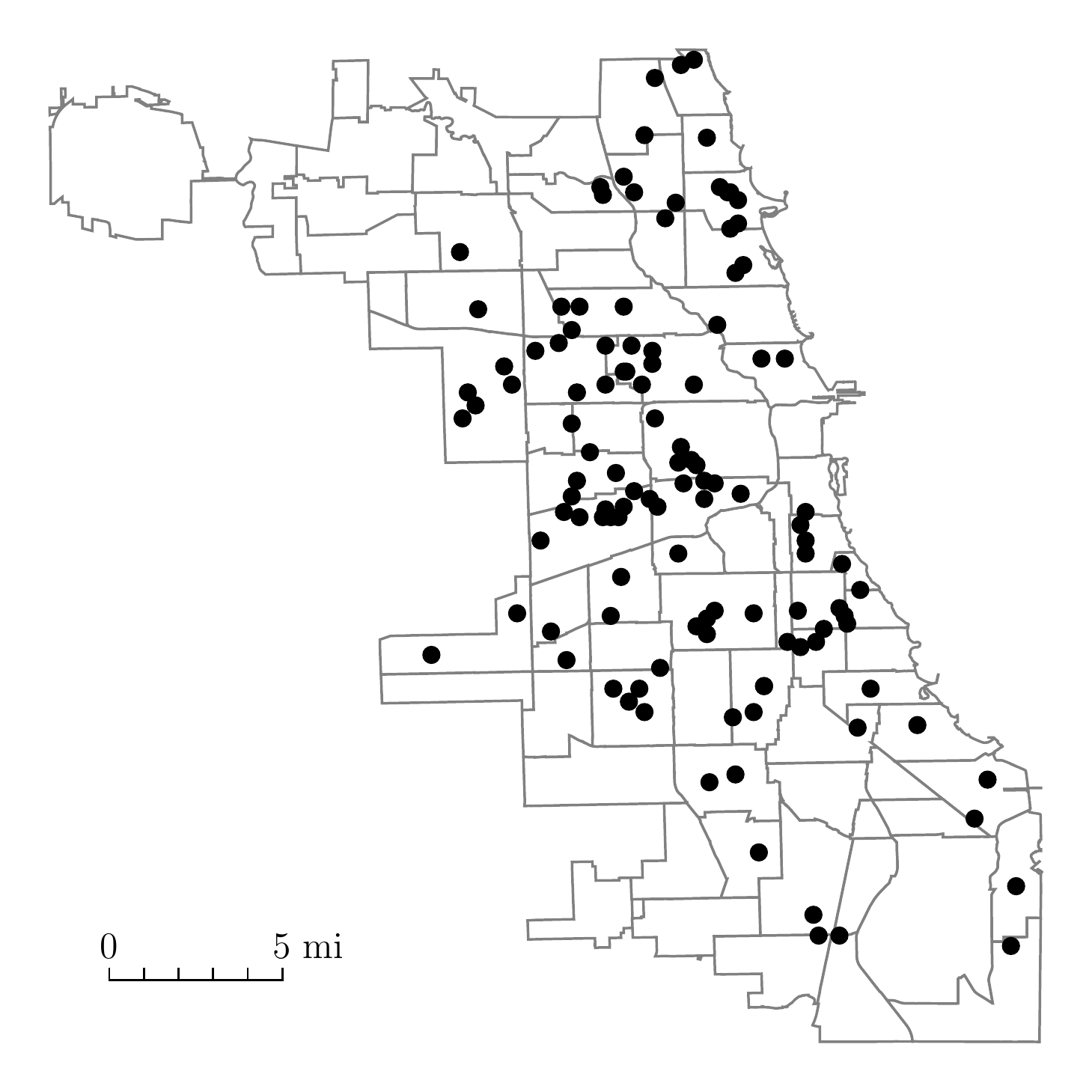}
	\caption{Map of the Chicago primary care community health clinic locations.}
	\label{fig:facility}
\end{figure}

A 2003 study by Luo and Wang \cite{luo2003b} applied various accessibility metrics to measure primary care accessibility throughout the Chicago region. They found that the distribution of accessibility levels can vary widely across even relatively small geographic areas of the city due to particular communities having an anomalously great or small amount of access relative to the surrounding areas. Specifically, their study identified three areas within the city with significantly better access than surrounding regions (The Loop, Lincolnwood-Skokie, and Elmhurst-Oak Brook) and two with significantly worse access (South Side and Midway), which agrees with the Chicago Department of Public Health findings. It should be noted that the South Side, in particular, is home to a disproportionately large number of citizens near or below the poverty level, many of whom lack health insurance and personal transportation. Figure \ref{fig:facility} shows the distribution of primary care community health clinics throughout the city \cite{primarycare}.

It is clear from the above studies that the city of Chicago is in need of more equitable access to primary care facilities. We will apply our modeling framework to examine changes that could potentially be made by the Chicago Transit Authority (CTA) in an effort to address this issue. Our goal is to find relatively simple, low-cost changes that could be put into practice with minimal impact on current service levels. Specifically, we consider altering the CTA bus network by reassigning vehicles between the existing bus lines, increasing the number of vehicles servicing certain lines while decreasing the number servicing others. We also consider the possibility of introducing express runs of existing routes.

%%%
\paragraph{Structure and Contributions of this Paper:}

In Section \ref{sec:model} we present the formulation of the optimization-based model, developed to fulfil the four modeling principles listed above, which constitutes the main contribution of this paper. We begin in Section \ref{subsec:overall} with its overall structure and then discuss its major components, including its underlying network representation in Section \ref{subsec:network}. In Section \ref{subsec:objective} we describe the model's objective function, which consists of a dynamic evaluation of the accessibility levels of the least-advantaged communities based on the current solution vector. In Section \ref{subsec:operator} we describe a set of model constraints which allow the planner to guarantee that operator and user costs remain within a specified margin of their initial levels. The structure of this model is unique in that its constraint functions are computationally expensive to evaluate while its objective function is not, and so in Section \ref{sec:algorithm} we present a hybrid tabu search/simulated annealing solution algorithm designed to minimize the number of computationally expensive evaluations needed.

The remainder of the paper is devoted to the computational trials used to test the performance of this model. In Section \ref{sec:trials} we explain the computational trials used as test cases, including the main CTA network case study in Section \ref{subsec:ctaproblem} and a set of artificial network tests in Section \ref{subsec:smallscaleproblem}. The results of these trials are presented in Section \ref{sec:results} and discussed in Section \ref{sec:discussion}. Appendix \ref{app:pseudocode} includes a more in-depth description of the solution algorithm from Section \ref{sec:algorithm}, Appendices \ref{app:data}--\ref{app:artificial} give detailed descriptions of the data and networks used for the computational trials explained in Section \ref{sec:trials}, and Appendix \ref{app:resultdata} includes the full data tables for these trials.

% Model formulation
\section{Model Formulation}
\label{sec:model}

In this section we describe our proposed model for addressing the design goals discussed above. We first show the overall structure of the model followed by detailed explanations of the individual components. Due to our choice of case study we will make specific reference to primary care access and a bus transit network, but these can easily be replaced with any other sort of objective and transit network that corresponds to the above assumptions.

%==============================================================================
\subsection{Overall Model Structure}
\label{subsec:overall}

In this section we give a high-level description of the overall structure of our proposed model. In the interest of keeping with our goal of considering only low-cost, immediate, easily-implemented tactical changes to the transit network, our design decisions consist of determining the daily average number of buses to assign to each line. A line's fleet size implicitly determines its average service frequency, and so this can be considered as a type of frequency-setting problem \cite{farahani2013}.

We begin by defining some variables that will be used throughout this paper. Let $L$ be the set of transit lines under consideration. For each line $l \in L$, let $y_l \in \mathbb{Z}$ be the number of vehicles that service line $l$, and let $\mathbf{y}$ be the vector of such variables arranged in a canonical order, which constitutes the decision vector for our design problem. Let $y_l^*$ be the number of vehicles initially assigned to line $l$, and $\mathbf{y^{\boldsymbol{*}}}$ be the initial fleet vector. Without loss of generality we assume that we are given constant upper and lower fleet size bounds $y_l^{\max}$ and $y_l^{\min}$, respectively, for each line $l$ (which can be made infinite if such bounds are not needed, or which can both be set to $y_l^*$ to force $y_l$ to remain constant). We assume that each line utilizes only a single type of vehicle, since a line that uses multiple vehicle types (e.g.\ a bus line serviced by both rigid-bodied buses and articulated buses) can be represented by multiple parallel lines in $L$. Let $Z$ be the set of all vehicle types utilized throughout the transit network. For each vehicle type $z \in Z$ let $L^z \subseteq L$ be the set of lines serviced by that vehicle, and for each line $l \in L$ let $z_l \in Z$ be the vehicle type used by line $l$.

Let $\mathbf{x}$ be the flow vector for user traffic through the transit network. User flows are defined implicitly as the output of a \textit{transit assignment model} which applies assumptions about user behavior to determine how the users, in aggregate, move through the system. The representation of the flow network and the transit assignment model will be discussed in Section \ref{subsec:network}. Because user decisions depend on the structure of the network as determined by the design decisions $\mathbf{y}$, we will denote the transit assignment function as $\mathrm{TransitAssignment}(\mathbf{y})$. Let $\mathbf{x^{\boldsymbol{*}}}$ be the user flow vector associated with the initial fleet vector.

In a typical public transit design model like those discussed above the objective is to minimize some combination of operator cost functions and user cost functions subject to design feasibility and budgetary constraints. We denote these two functions as $\mathrm{OperatorCost}(\mathbf{y},\mathbf{x})$ and $\mathrm{UserCost}(\mathbf{x})$, respectively. It may be assumed that operator costs depend directly on both our design decision vector $\mathbf{y}$ and the user flow vector $\mathbf{x}$, but the user costs generally depend directly only on the user flows (though both depend implicitly on the design decisions). Possible choices for these functions will be discussed in Section \ref{subsec:operator}, but in general operator costs represent costs associated with running the transit network while user costs represent how difficult it is for users to travel through the network.

Finally, let $\mathrm{Access}(\mathbf{y})$ be a function describing the social access objective. Possible choices for this objective will be discussed in Section \ref{subsec:objective}, but in general it should represent how easy it is for users to access some service within the network, weighted to focus on the most disadvantaged users. It depends only on the network design $\mathbf{y}$ and not the user flows $\mathbf{x}$ because, in line with modeling principle 4, the user flow volume is assumed to be small and infrequent enough that it does not compete for transit capacity in the same way that day-to-day flow volume does.

The three objectives $\mathrm{OperatorCost}(\mathbf{y},\mathbf{x})$, $\mathrm{UserCost}(\mathbf{x})$, and $\mathrm{Access}(\mathbf{y})$ cannot, in general, all be optimized simultaneously, which makes this a multiobjective program. We take the approach of optimizing only the access objective while including the operator and user costs as constraints. This allows our model to explicitly attempt to make improvements to the social access objective while still guaranteeing a specified service level. Our proposed model takes the form of the following nonlinear program, which we will refer to as the \textit{social access maximization problem} (SAMP):
\begin{alignat}{2}
	\label{eqn:nlp1objective} \max_\mathbf{y} \quad& \mathrm{Access}(\mathbf{y}) \\
	\label{eqn:nlp1operator} \mathrm{s.t.} \quad& \mathrm{OperatorCost}(\mathbf{y},\mathbf{x}) \le B_\text{operator} \\
	\label{eqn:nlp1user} & \mathrm{UserCost}(\mathbf{x}) \le B_\text{user} \\
	\label{eqn:nlp1assignment} & \mathbf{x} = \mathrm{TransitAssignment}(\mathbf{y}) \\
	\label{eqn:nlp1vehicle} & \sum_{l \in L^z} y_l \le \sum_{l \in L^z} y_l^* &\qquad& \forall z \in Z \\
	\label{eqn:nlp1bounds} & y_l^{\min} \le y_l \le y_l^{\max} && \forall l \in L \\
	\label{eqn:nlp1integer} & y_l \in \mathbb{Z} && \forall l \in L
\end{alignat}

To explain, the overall objective (\ref{eqn:nlp1objective}) is to maximize the access objective. Constraints (\ref{eqn:nlp1operator}) and (\ref{eqn:nlp1user}) give explicit upper bounds of $B_\text{operator}$ and $B_\text{user}$ to the operator and user costs, respectively. These bounds allow the transit planner to specify precisely how much of an increase is acceptable in exchange for maximizing the social access objective, which can be chosen based on their requirements. Constraints (\ref{eqn:nlp1assignment}) implicitly define the user flow vector $\mathbf{x}$ as a result of the decision vector's effect on the transit assignment function. Finally constraints (\ref{eqn:nlp1vehicle}) ensure that no new vehicles of any type are added, constraints (\ref{eqn:nlp1bounds}) enforce fleet size bounds for each line, and constraints (\ref{eqn:nlp1integer}) enforce integrality of fleet sizes. Note that, in general, transit assignment functions are computationally expensive to evaluate, and so determining the user flows and thus the feasibility of constraints (\ref{eqn:nlp1operator}) and (\ref{eqn:nlp1user}) for a given solution vector $\mathbf{y}$ is expected to be computationally expensive. As a result it is far faster to evaluate the objective of a solution to the SAMP than its feasibility, and this unusual structure dictates the form of the solution algorithm discussed in Section \ref{sec:algorithm}.

As the SAMP is meant to represent a general modeling framework for handling social access-oriented public transit design problems, the exact definitions of the functions $\mathrm{Access}$, $\mathrm{OperatorCost}$, $\mathrm{UserCost}$, and $\mathrm{TransitAssignment}$ should be chosen based on planner's specific goals. The remainder of this section will be spent describing some possible choices for these functions and explain the modeling decisions made for our case study of Chicago primary care access.

%==============================================================================
\subsection{Transit Assignment Model and Transit Network Representation}
\label{subsec:network}

The mathematical representation of the transit system takes the form of a network with nodes to represent locations (like transit stops, population centers, and primary care facilities) and arcs to represent traversable routes between those locations (like bus trips, train trips, and sidewalk travel), but the specific structure of this network depends on the choice of transit assignment model, which may require the inclusion of auxiliary nodes and arcs to represent different parts of the user decision-making process. In this section we describe our choice of transit assignment model and the network representation that it requires.

%%%
\paragraph{Assignment Model:} The transit assignment model, $\mathrm{TransitAssignment}(\mathbf{y})$, applies a set of basic assumptions regarding user behavior to simulate how public transit users route themselves through the transit network. The flow produced by such a model is called the \textit{user-optimal flow}, which in general leads to worse objective values than would be possible if $\mathbf{x}$ were directly controllable. There is an extremely wide variety of transit assignment models used throughout transportation design literature, differing primarily in their assumptions regarding user behavior and how they account for user interactions in the presence of congestion  \cite{fu2012}. For the purposes of our study we have chosen to use a nonlinear model by Spiess and Florian \cite{spiess1989}. This model works by using a conical congestion function to increase the costs of arcs as they become more crowded, serving to discourage flows from exceeding line capacities. The output is the user-optimal flow vector which corresponds to all users choosing a strategy for for which their expected travel time cannot be improved with a unilateral change. This particular model was chosen because it is well-established and commonly used in transit forecasting systems (such as EMME/2) as well as being relatively simple, which is desirable since the solution algorithm presented in Section \ref{sec:algorithm} potentially requires many evaluations of the operator cost and user cost functions during each iteration.

%%%
\paragraph{Network Structure:} Our network representation of the public transit system was chosen to meet the requirements of this assignment model. We use directed graph $G = (V,E)$ to represent the public transit network. The node set $V$ describes locations of interest within the city, including origins, destinations, and public transit stops. The arc set $E$ represents links between these locations via different modes of transportation, including walking, bus travel, and train travel. We consider only public transit and foot travel, not automobile traffic.

We partition the node set into subsets. Let $V_\text{stop} \subset V$ be the set of \textit{stop nodes}, which represent the physical locations of public transit stops. Let $V_\text{board} \subset V$ be the set of \textit{boarding nodes}, which are used to distinguish between different transit lines at each stop. For each stop node $i \in V_\text{stop}$ there is a boarding node $j$ for all lines $l \in L$ that utilize stop $i$. Let $V^+ \subset V$ be the set of \textit{origin nodes}, which represent locations the locations where public transit trips begin. Finally, let $V^- \subset V$ be the set of \textit{destination nodes}, which represent locations the locations where public transit trips end. Note that these origin-destination (OD) pairs represent the high-volume, high-frequency day-to-day trips that the assignment model considers, not the social access-related objectives.

We similarly partition the arc set. Let $E_\text{line} \subset E$ be the set of \textit{line arcs}, which represent in-vehicle travel for a transit line. Let $E_\text{board} \subset E$ be the set of \textit{boarding arcs}, which represent the decision to board a particular line at a particular stop. Let $E_\text{alight} \subset E$ be the set of \textit{alighting arcs}, which represent the decision to alight a particular line at a particular stop. Finally, let $E_\text{walk} \subset E$ be the set of \textit{walking arcs}, which represent travel along sidewalks or other walking paths.

\begin{figure}[h]
	\centering
	\includegraphics[width=0.6\textwidth]{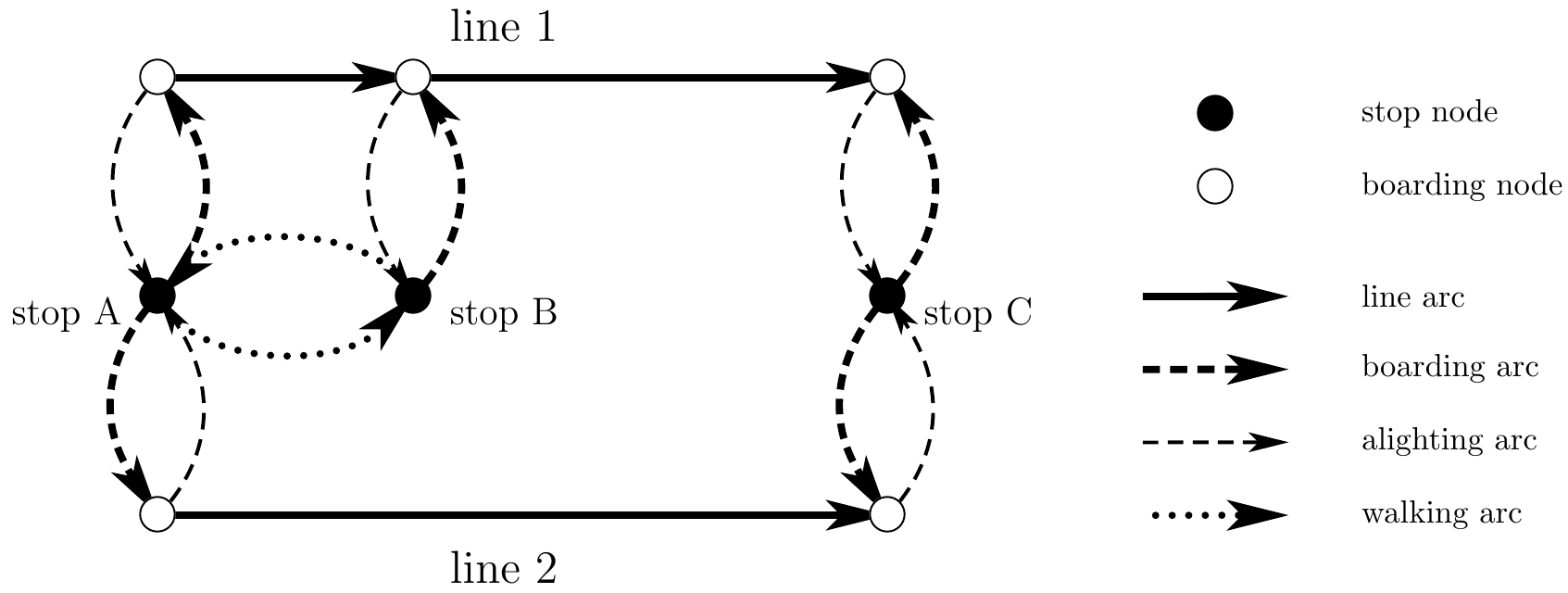}
	\caption{Illustration of the structure of the transit network directed graph. This example shows a sequence of three stops (A, B, and C) serviced by two transit lines (1 and 2). Line 1 stops at all three locations while line 2 stops at only A and C. Stops A and B are within walking distance of each other.}
	\label{fig:netstructure}
\end{figure}

For each line $l \in L$, let $G^l = (V^l,E^l)$ be the subgraph of $G$ containing only nodes and arcs associated with that line. Let any of the above partitions with a superscript of $l$ represent the subset of that partition corresponding to line $l$ (so for example $V_\text{stop}^l$ is the set of all stop nodes traversed by line $l$). Note that origin nodes, destination nodes, and walking arcs have no associated lines. For every stop node $i \in V_\text{stop}^l$ traversed by line $l$, there is a corresponding boarding node $i_l \in V_\text{board}^l$. These two nodes are joined by a boarding arc $(i,i_l) \in E_\text{board}^l$ and an alighting arc $(i_l,i) \in E_\text{alight}^l$. If stop node $j \in V_\text{stop}^l$ is the next stop on line $l$ after stop $i$, and if its corresponding boarding node is $j_l \in V_\text{board}^l$, then there is a line arc $(i_l,j_l) \in E_\text{line}^l$ joining the consecutive boarding nodes. Walking arcs are not line-specific and may join any pairs of stop nodes, origin nodes, and destination nodes. Figure \ref{fig:netstructure} shows an example of this network structure.

In order to incorporate the social access objective into this network, we also define two more types of node. Let $\widetilde{V}^+ \subset V$ be the set of \textit{community nodes}, which represent locations the communities of concern in the social access objective. Let $\widetilde{V}^- \subset V$ be the set of \textit{facility nodes}, which represent locations of the service of concern in the social access objective. Whereas the objective of a typical public transit study is to facilitate travel between origin nodes $V^+$ and destination nodes $V^-$, the social access objective is to facilitate access of nodes $\widetilde{V}^+$ to nodes in $\widetilde{V}^-$. Like origin and destination nodes, walking arcs connect community and facility nodes to the rest of the network.

%%%
\paragraph{Frequency-Based Attributes:} To relate the structure of network $G$ to the SAMP we use a frequency based approach, a common choice for many tactical UTNDP problems \cite{guihaire2008,farahani2013} and transit assignment models \cite{kurauchi2003,marcotte2004,cepeda2006,schmocker2011,fu2012,schmocker2013}, including our chosen model. Under this paradigm all traffic demands and flows are averaged over a specific time period during which it is assumed that travel demand and traffic levels remain constant. The frequency-based approach does not take vehicle timetables into consideration, instead considering only the expected arrival frequency of each line at each stop which in turn determines a user's expected waiting time. We will use $\tau$ to represent the daily time horizon for use in the frequency-based model.

For each OD pair $(s,t) \in V^+ \times V^-$, let $D_{st}$ be the total day-to-day transit demand over the specified time period. The $|V^+| \times |V^-|$ matrix for which element $(s,t)$ is $D_{st}$ is known as the \textit{OD matrix}. For each arc $ij \in E$, let $x_{ij}$ be the total day-to-day flow through that link, and let $\mathbf{x}$ be the vector of such flow values (ordered in a canonical way).

The decision variables $y_l$ of the SAMP represent the average number of vehicles servicing each line $l \in L$ on a day-to-day basis. We assume that these tactical decisions are minor enough not to disturb any overall traffic patterns, and so each line also has a constant expected time of $\tau^l$ (in minutes) to complete a full circuit (including layovers). These decision variables affect the output of the transit assignment model in two ways: First, the number of vehicles servicing a line dictates that line's average service frequency $f_l$ (in units of arrivals/minute) which, under the assumption of exponentially distributed arrival times, is determined by the inverse proportionality $f_l = y_l/\tau^l$ \cite{spiess1989}. Second, the number of vehicles servicing a line dictates that line's total capacity $u_l$ (in units of total passengers over the given time horizon), which for a line with a vehicle seating capacity of $b_l$ is determined by $u_l = \tau b_l f_l$. Within the assignment model line frequency is used to determine expected waiting time while line capacity is used to determine penalty costs for overcrowding, both of which contribute to users' path selections. It should be noted that any line which is assigned no vehicles, $y_l = 0$, has an effective capacity of zero and an effectively infinite waiting time, causing it to be chosen by no users.

For the purposes of the social access objective we assume that the volume of traffic is so small that it makes up an insignificant fraction of the day-to-day travel volume, and so it is ignored within the assignment model. For this reason it can be assumed not to depend on the capacities of the lines, but it may still depend on the frequencies, since they determine expected waiting times and thus expected trip length. In addition, the introduction of new lines such as the express runs to be discussed in Section \ref{subsec:ctaproblem} can drastically affect both day-to-day traffic and social access traffic by making new paths available.

%==============================================================================
\subsection{Objective Function}
\label{subsec:objective}

There are many reasonable choices for the social access objective function $\mathrm{Access}(\mathbf{y})$, but in general it is meant to describe the ease of traveling from a community to a facility using the network resulting from design decisions $\mathbf{y}$, and to do so in a way that places focus on equity of access. Our proposed objective is based on maximizing primary care accessibility metrics for the set of communities that currently has the lowest level of access. As referenced in Section \ref{sec:intro}, a large amount of research has gone into the development of accessibility metrics, particularly with regard to primary care access. For example, the \textit{two-step floating catchment area} (2SFCA) metric was first described in Luo and Wang 2003 \cite{luo2003a} and has been used in a wide variety of studies since then \cite{luo2003b,wang2005,luo2009,mcgrail2009,mcgrail2014,neutens2015}.

Let $A_i(\mathbf{y})$ be a primary care accessibility metric of community $i \in \widetilde{V}^+$. Our specific choice of $A_i(\mathbf{y})$ will be described below, but in general accessibility metrics depend on travel times within the network, and so are implicitly functions of the decision vector $\mathbf{y}$. Then we define our objective function as the sum of the $\mathcal{K}$ minimum accessibility metrics,
\begin{align}
	\label{eqn:accesscost} \mathrm{Access}(\mathbf{y}) &:= \sum_{\substack{\text{$\mathcal{K}$ minimum elements} \\ \text{of $\{ A_i(\mathbf{y}) \, | \, i \in \widetilde{V}^+\}$}}} A_i(\mathbf{y})
\end{align}

where $\mathcal{K}$ is a model parameter. Note that as the solution algorithm proceeds and the solution vector $\mathbf{y}$ is changed, the accessibility metrics $A_i(\mathbf{y})$ and their relative rankings may also change. This objective refers to the $\mathcal{K}$ \textit{current} least accessibility metrics (for the current solution $\mathbf{y}$).

This objective generalizes the ideas of maximizing only the current minimum accessibility metric (which is the special case of $\mathcal{K}=1$), and of maximizing the total (or equivalently the average) of all accessibility metrics (which is the special case of $\mathcal{K}=|\widetilde{V}^+|$). Attempting to optimize for only the current least-accessible community leads to potential problems if that community cannot be significantly helped through changes in the transit network, since it results in the solution process focusing all of its decisions on attempting to make minor improvements in one location while ignoring decisions that might help far more communities. On the other hand, attempting to optimize the total accessibility level of all communities could result in decisions which improve the accessibility levels of areas that already have high accessibility while worsening the accessibility levels of disadvantaged areas, as long as it improves the overall total. The computational trials for our Chicago case study, as well as for artificial test networks, displayed both of these problems for extreme values of $\mathcal{K}$ (see Section \ref{sec:results}).

Our specific choice of accessibility metric $A_i(\mathbf{y})$ is a \textit{gravity-based metric}. Gravity-based metrics take into account how many facilities are near a community, however rather than a binary count of whether or not the travel time falls below a cutoff, each facility is weighed by a negative power of the pairwise travel time, causing the supposed access level to a facility to decay over distance like a gravitational force. This was chosen for use in the SAMP objective because, unlike the catchment area-based metrics, its continuous dependence on travel times generally gives it a nonzero gradient, which helps to give direction to the the solution algorithm's local search (see Section \ref{sec:algorithm}).

The specific metric used in our study is a competition-based model developed in Weibull 1976 \cite{weibull1976} for describing differential access to employment opportunities. It has also been used in the study of health care access, for example in Joseph and Bantock 1982 \cite{joseph1982} and in Luo and Wang 2003 \cite{luo2003b}. It is defined as
\begin{alignat}{2}
	\label{eqn:gravitymetric} A_i(\mathbf{y}) &:= \sum_{j \in \widetilde{V}^-} \frac{S_j d_{ij}^{-\beta}(\mathbf{y})}{F_j(\mathbf{y})} &\qquad& \forall i \in \widetilde{V}^+
\end{alignat}

where $S_j$ is the quality of facility $j \in \widetilde{V}^-$ (usually taken as its patient capacity), $d_{ij}(\mathbf{y})$ is the travel time from $i$ to $j$ (calculated as graph distance in $G$, which depends on the network design decisions $\mathbf{y}$), $\beta > 0$ is a gravitational decay model parameter, and $F_j(\mathbf{y})$ is a competition-related facility metric defined by
\begin{alignat}{2}
	\label{eqn:gravitycompetition} F_j(\mathbf{y}) &:= \sum_{k \in \widetilde{V}^+} P_k d_{kj}^{-\beta}(\mathbf{y}) &\qquad& \forall j \in \widetilde{V}^-
\end{alignat}

where $P_k$ is the population of community $k$ seeking service at facility $j$. $F_j(\mathbf{y})$ can be interpreted as a measure of how crowded facility $j$ is, and is greater if the facility is close to many populous communities. The accessibility metric $A_i(\mathbf{y})$ rewards a community for being close to many facilities with high quality and low overcrowding. The metric, itself, has no direct interpretation in isolation and is meant for use in comparing communities across space or time, as will be shown in Section \ref{sec:results}.

%==============================================================================
\subsection{Operator Costs and User Costs}
\label{subsec:operator}

The operator and user cost functions represent typical objectives from standard UTNDP studies. Operator costs generally represent monetary costs incurred by the public transit authority, and may include vehicle operation costs, vehicle acquisition costs, and fare revenue. User costs generally represent time costs incurred by users of the transit system, and may include travel time (which may be weighted by mode), waiting time, number of transfers, and measures of discomfort.

Because this study involves only the reassignment of existing vehicles within the public transit network and not the purchase of new vehicles, and because operator costs consist largely of budgetary concerns, the operator cost function is unlikely to change as a result of our decisions. For this reason we will exclude the operator cost constraint (\ref{eqn:nlp1operator}) from our model for the remainder of this paper, however it would still be appropriate to include for a different set of design decisions (for example, including vehicle purchase or line construction). Our user cost function takes the form of the weighted sum
\begin{align}
	\label{eqn:usercost} \mathrm{UserCost}(\mathbf{x},\omega) &:= \theta_1 \sum_{ij \in E_\text{line}} c_{ij} x_{ij} + \theta_2 \sum_{ij \in E_\text{walk}} c_{ij} x_{ij} + \theta_3 \omega
\end{align}

To explain, although the form of the user cost function mentioned in Section \ref{subsec:overall} included only a flow vector argument $\mathbf{x}$, our chosen assignment model produces both a flow vector and a total waiting time scalar $\omega$ which represents the total time spent by all users waiting to board \cite{spiess1989}. The above user cost function is then a weighted sum of three types of travel time: total in-vehicle travel time, total walking travel time, and total waiting time, with weights of $\theta_1$, $\theta_2$, and $\theta_3$, respectively, which can be used to penalize particular types of travel time more harshly than others. The user cost bound constraint (\ref{eqn:nlp1user}) in the SAMP prohibits $\mathrm{UserCost}(\mathbf{x},\omega)$ from increasing beyond a constant upper bound $B_\text{user}$. We use a constant bound of
\begin{align}
	\label{eqn:usercostbound} \mathrm{UserCost}(\mathbf{x},\omega) \le (1+\epsilon) \, \mathrm{UserCost}(\mathbf{x^{\boldsymbol{*}}},\omega^*)
\end{align}

where $\epsilon \ge 0$ is a model parameter which specifies the allowed margin of increase in user costs, $\mathbf{x^{\boldsymbol{*}}}$ is the initial user flow vector, and $\omega^*$ is the initial waiting time scalar. This bound corresponds to the choice of $B_\text{user} = (1+\epsilon) \, \mathrm{UserCost}(\mathbf{x^{\boldsymbol{*}}},\omega^*)$ and ensures that the user cost cannot undergo a relative increase of more than $\epsilon$ above its initial level. A larger value of $\epsilon$ results in a less restrictive bound and thus a larger feasible set and potentially a better optimal value, but at the cost of potentially increasing the overall user costs. In this way $\epsilon$ allows for different choices in the tradeoffs between improving social access and maintaining day-to-day service levels.

% Solution algorithm
\section{Hybrid Tabu Search/Simulated Annealing Solution Algorithm}
\label{sec:algorithm}

In this section we give a brief description of our proposed solution method for the SAMP model defined in Section \ref{sec:model}. As mentioned above, we wish to allow for a variety of choices in the definitions of $\mathrm{Access}(\mathbf{y})$, $\mathrm{OperatorCost}(\mathbf{y},\mathbf{x})$, $\mathrm{UserCost}(\mathbf{x})$, and $\mathrm{TransitAssignment}(\mathbf{y})$, and so we will treat all of these functions as black boxes, assuming that each can be solved through some subroutine and making use only of its outputs for a given set of inputs.

Versions of the UTNDP are well-known as computationally difficult problems for a variety of reasons, including their large number of integer variables and the bilevel structure which arises from the need to consider user behavior as an implicit function of the design decisions. As a result they are in general nonlinear and nonconvex, and most studies on the subject resort to a heuristic or metaheuristic solution method to obtain high-quality solutions within a reasonable time frame \cite{bielli2002,chakroborty2002,bagloee2011,cipriani2012,farahani2013}. Transit assignment models are in general computationally expensive to solve \cite{guihaire2008}, and so it is reasonable to assume that $\mathrm{TransitAssignment}(\mathbf{y})$ and anything that relies on its outputs (including $\mathrm{OperatorCost}(\mathbf{y},\mathbf{x})$ and $\mathrm{UserCost}(\mathbf{x})$) are significantly more expensive to evaluate than the other major components of the model. This gives the SAMP the unusual structural property of having constraints (\ref{eqn:nlp1operator})--(\ref{eqn:nlp1user}) that are significantly more expensive to evaluate than the objective (\ref{eqn:nlp1objective}). Because of this, our proposed metaheuristic approach is designed to minimize the number of required constraint function evaluations.

We propose a hybrid tabu search/simulated annealing solution algorithm for the SAMP, due to the fact that there is a natural choice of initial feasible solution (the initial fleet vector, $\mathbf{y^{\boldsymbol{*}}}$) as well as a natural definition of local moves, consisting of adding, dropping, and swapping individual vehicles between compatible lines. Tabu search (TS) and simulated annealing (SA) are widely-used combinatorial optimization metaheuristics. They are both based on local search, wherein each solution has a defined \textit{neighborhood} of solutions which differ by only one local move, and the solution algorithm consists of iteratively conducting neighborhood searches and making local moves to improve the objective value. Pure local search carries the risk of becoming trapped at a local optimum while much better solutions exist, and so TS and SA both include rules for occasionally making suboptimal moves in order to explore more of the feasible set and arrive at a better solution.

SA accomplishes this by introducing randomness \cite[pp.~512--514]{linopt}, having a chance to make suboptimal moves according to a parameter $T$ called \textit{temperature} which decreases as the algorithm runs. As $T$ decreases suboptimal moves become less likely, and so late in the run time it is more likely to remain near local optima. On the other hand, TS accomplishes this by using memory structures to avoid repeatedly cycling through the same solutions \cite{glover2013}, using a combination of \textit{short-term memory} (STM) to avoid reversing a move too soon after it is made, and \textit{long-term memory} (LTM) to remember profitable past solutions that can be returned to if the algorithm becomes trapped. The STM generally consists of a list of \textit{tabu rules} which exist for a set \textit{tenure} before being dropped from the tabu list. Hybrids of the two algorithms have been shown in a number of recent combinatorial optimization studies to perform better than pure TS or SA alone \cite{arikan2012,alhroob2014,lenin2016,shafahi2018}.

The general structure of our approach is based primarily on that of Alhroob et al.~2014 \cite{alhroob2014}. It follows the TS structure discussed above, but uses an SA acceptance probability when considering whether to make a suboptimal move. Tabu rules stored in STM prevent undoing recent additions or subtractions from each fleet, and attractive solutions stored in LTM consist of a combination of the second-best solutions from neighborhood searches which were not chosen, and of suboptimal moves that were previously denied by the SA criterion. Tabu tenures are adjusted adaptively to increase when more diversification is required and to decrease when more intensification is required, measured by using two separate nonimprovement move counters: An inner counter increments along with the tabu tenures whenever no improving non-tabu move exists, resets whenever a suboptimal move is made by passing the SA criterion, and backtracks to a random attractive solution from the LTM when it reaches a set maximum. An outer counter increments alongside the inner counter as well as whenever a an attractive solution is backtracked to, resets whenever an improving local move is made, and prompts the tabu tenures to reset when it reaches a set maximum.

In order to take advantage of the fact that the social access objective can be evaluated relatively quickly while the operator and user cost constraints cannot, our algorithm uses a two-pass approach: The first pass consists of considering every possible ADD and DROP move that satisfies the trivial design constraints (\ref{eqn:nlp1vehicle})--(\ref{eqn:nlp1integer}), evaluating its objective value, considering SWAP moves made up of the most profitable ADD and DROP moves combined, and evaluating their objective values. The second pass then consists of evaluating the much more computationally expensive user flow-dependent feasibility constraints (\ref{eqn:nlp1operator})--(\ref{eqn:nlp1user}) of the neighbors in descending order of objective value, halting as soon as two feasible solutions are found, which due to the first pass are guaranteed to be the two local moves with the maximum objective value. At all times tabu rules are enforced to limit the neighborhood, unless a new global optimum is found.

The following is a brief summary of our proposed solution algorithm. The full pseudocode can be found in Appendix \ref{app:pseudocode}.

\begin{enumerate}[1.]
	\small
	\item Initialize solution vector $\mathbf{y}^{(0)}$ as the current real-world fleet size vector $\mathbf{y^{\boldsymbol{*}}}$, and set it as the best known solution. Initialize TS and SA memory structures, including the SA temperature $T$, tabu tenures $t$, empty STM and LTM sets, and zero nonimprovement counters.
	\item Repeat the following for a predetermined number of iterations $k=0,1,\dots,K$:
	\begin{enumerate}[a.]
		\item (ADD/DROP Neighborhood Search, First Pass) Consider every possible line $l \in L$ in a random order. If adding one vehicle to line $l$ is not tabu (or if adding to $l$ offers an improved best solution), and if it does not violate the vehicle limit bounds (\ref{eqn:nlp1vehicle})--(\ref{eqn:nlp1integer}), save $l$ in a tentative ADD move list. Repeat an analogous process to gather a tentative DROP move list.
		\item (ADD/DROP Neighborhood Search, Second Pass) Consider all lines $l$ collected in the tentative ADD move list, in descending order of objective value (i.e.\ from best to worst). For each line, evaluate user and operator constraints (\ref{eqn:nlp1operator})--(\ref{eqn:nlp1user}). If the solution is feasible, add it to a final ADD move list. Stop as soon as a desired number of ADD moves has been found. Repeat an analogous process to gather a final DROP move list.
		\item (SWAP Neighborhood Search) Consider all pairs of lines $(m,l)$ where $m$ is a finalized DROP move and $l$ is a finalized ADD move, in descending order of the total objective value of the two individual moves. If swapping a vehicle from line $m$ to line $l$ does not violate the user and operator constraints (\ref{eqn:nlp1operator})--(\ref{eqn:nlp1user}), add the move to a final SWAP move list. Stop as soon as a desired number of SWAP moves has been found.
		\item Among the finalized ADD, DROP, and SWAP move lists, select the two solutions with the greatest objective values. Let $\mathbf{\tilde{y}}^{(k,1)}$ be the best and $\mathbf{\tilde{y}}^{(k,2)}$ be the second-best.
		\item If $\mathbf{\tilde{y}}^{(k,1)}$ is an improvement over the current solution $\mathbf{y}^{(k)}$, then move to the improving neighbor by defining $\mathbf{y}^{(k+1)}$ as $\mathbf{\tilde{y}}^{(k,1)}$, reset the inner nonimprovement counter and tabu tenures, and add tabu rules to STM to prevent undoing any ADD or DROP that was just conducted.
		\item Otherwise, increment both nonimprovement counters and then apply do the following:
		\begin{enumerate}[i.]
			\item (SA Acceptance Criterion) With probability $\exp\{-[ \mathrm{Access}(\mathbf{y}^{(k)}) - \mathrm{Access}(\mathbf{\tilde{y}}^{(k,1)})]/T\}$, move to the best neighbor by defining $\mathbf{y}^{(k+1)}$ as $\mathbf{\tilde{y}}^{(k,1)}$, increment the tabu tenures, reset the inner nonimprovement counter, define tabu rules in STM to prevent undoing any ADD or DROP that was just conducted, and save the second-best neighbor $\mathbf{\tilde{y}}^{(k,2)}$ in LTM.
			\item Otherwise, save the best neighbor $\mathbf{\tilde{y}}^{(k,1)}$ in LTM.
		\end{enumerate}
		\item (Evaluate Nonimprovement Counters) If the inner nonimprovement counter has met its maximum, then reset it, increment the outer nonimprovement counter, and move to a random solution from LTM. If the outer nonimprovement counter has met its maximum, reset the tabu tenures.
		\item (Upkeep) Drop moves from STM whose tabu tenures have elapsed. Apply a cooling schedule to the temperature $T$.
	\end{enumerate}
	\item Return the best known solution and its objective value.
\end{enumerate}

The source code for our C++ implementation can be viewed online \cite{p2codemaster}. The algorithm includes a large number of tuning parameters which were adjusted using trial-and-error during the computational trials described below. For the purposes of the Chicago case study, the SA criterion used an initial temperature of $10^{16}$ and a cooling schedule consisting of multiplying the temperature by a factor of $0.999$. We allowed a maximum of 40 solutions in LTM, after which attractive solutions were randomly dropped. We evaluated $50$ ADD and DROP moves during each first pass and kept $2$ of each during each second pass. The initial tabu tenures were $6.0$ and increased by a factor of $1.15$ whenever incremented. Nonimprovement counter cutoffs of $20$ (inner) and $10$ (outer) were used.

% Computational trials
\section{Computational Trials}
\label{sec:trials}

In this section we describe the computational trials used for evaluating the performance and results of the SAMP. In addition to the Chicago primary care access case study described in Section \ref{sec:intro}, we also consider a small-scale artificial network in order to generate additional computational results for use in sensitivity analysis. The results of these tests are given in Section \ref{sec:results}.

%==============================================================================
\subsection{CTA Design Problem}
\label{subsec:ctaproblem}

In this section we briefly explain the setup of the problem and some of the tests that were run using the model. An in-depth explanation of the data and network setup is given in Appendix \ref{app:data}. The main design problem of interest for this study is that of making minor alterations to the CTA network in order to improve equity of access to primary care facilities. To that end an underlying network $G = (V,E)$ was constructed to represent the CTA public transit network. The set of stop nodes $V_\text{stop}$ corresponded to the set of public transit stops within the city and were also used as the origins and destinations for day-to-day transit volume. The OD matrix was estimated from stop-level boarding data. To reduce the size of the network we clustered the more than $11,000$ defined CTA stops by identifying nearby locations with each other, resulting in a set of $|V_\text{stop}| = 997$ stops.

Community nodes $\widetilde{V}^+$ were taken as population centers within the city, defined using 2010 census data clustered into population-weighted centroids for each of the 77 Chicago community areas, resulting in $|\widetilde{V}^+| = 77$ community nodes. Facility nodes $\widetilde{V}^-$ were taken as primary care facilities, obtained from a list of primary care community health clinics provided by the City of Chicago, resulting in $|\widetilde{V}^-| = 119$ facility nodes (see Figure \ref{fig:facility}). Each facility $j \in \widetilde{V}^-$ was given a constant quality metric of $S_j \equiv 1$ in (\ref{eqn:gravitymetric}). Walking arcs were generated to connect locations within $0.75$ miles of each other, with travel times calculated under the assumption of 4 ft/sec walking speed. The final network $G$ included a total of 3852 nodes and 17,522 arcs.

The public transit lines were defined using data available from CTA General Transit Feed Specification (GTFS) files. The network includes 126 bus lines and 8 train lines. Schedule data were used to estimate the current number of vehicles $y_l^*$ as well as the average circuit time $\tau^l$ for each line $l \in L$. The unit time horizon $\tau$ was taken as the full 24-hour (1440 minute) period. For simplicity we assumed that all bus lines used a single type of bus (the New Flyer D40LF, with $b_l = 39$ seats per bus) and that all train lines used a single type of train (a 6 car train of Bombardier Transportation 5000 series cars, with $b_l = 228$ seats per train). Throughout our computational trials we considered only changes to the bus lines, with the train lines being kept constant. The bus fleet lower bounds $y_l^{\min}$ were taken as the minimum number of vehicles required to achieve a line frequency $f_l$ of at least one arrival per 30 minutes. No upper bounds were used.

In order to evaluate whether the introduction of express lines can positively affect primary care access levels, we also defined potential express runs of existing bus lines. An express run follows the same route as the original line but stops at only a subset of the original stops, reducing its circuit completion time $\tau^l$ and thus increasing its effective frequency $f_l$ as well as some of its pairwise travel times. While it is possible to build express line design decisions into the SAMP by introducing binary decisions for which stops to include and exclude, we instead elected to generate a small set of heuristically designed express lines (see Appendix \ref{app:expressdata}). These lines were built into $G$ alongside the current CTA lines but defined to have initial fleet sizes of $y_l^* \equiv 0$, making them unavailable for use until adding at least one vehicle. The bus reassignment decisions $\mathbf{y}$ then allow the planner to transfer buses from existing lines to the new express lines if doing so would improve the access objective. A total of 65 express runs of existing lines were generated, each of which included a subset of approximately 20\% of the original line's stops chosen to connect areas of high access to areas of low access.

The solution algorithm discussed in Section \ref{sec:algorithm} was used to solve the instances of the SAMP discussed above, running for 500 search iterations followed by an exhaustive neighborhood search of the best known solution to ensure local optimality.
Throughout all tests the access function (\ref{eqn:accesscost}) was evaluated using a gravitational decay parameter of $\beta = 1.0$ and a community count of $\mathcal{K} = 8$ (which is approximately $10.39\%$ of the 77 community areas), meaning that the objective was to maximize the 8 lowest accessibility metrics. The user cost function (\ref{eqn:usercost}) used equal mode weights of $\theta_1 = \theta_2 = \theta_3 = 1$ and a maximum relative increase of $\epsilon = 0.01$, with constraint (\ref{eqn:usercostbound}) preventing total user travel time from increasing by more than $1\%$ above initial levels. The results of these trials are discussed in Section \ref{subsec:ctaresults}. A few modifications of the main trial were also run, including: a version which ignored the user cost constraint (\ref{eqn:nlp1user}), a version for which the only allowed moves were to transfer vehicles between lines and their express runs, and a version for which the express lines were excluded. The results of these supplemental trials are discussed in Section \ref{subsec:ctasupplemental}.

%==============================================================================
\subsection{Artificial Network Tests}
\label{subsec:smallscaleproblem}

Part of our evaluation of the SAMP model involved conducting sensitivity analysis to determine whether its results depend strongly on the particular choice of modeling parameters $\beta$, $\mathcal{K}$, and $\epsilon$. The full CTA problem instance is too computationally expensive to solve repeatedly under varying circumstances, and so in order to generate more data we defined a small artificial network whose SAMP instance could be solved much more quickly.

A full description of the artificial test network is given in Appendix \ref{app:artificial}. It represents a $0.9 \times 1.5$ mile city arranged in a square grid of blocks on a $7 \times 11$ street lattice, and contains a total of 527 nodes (including 204 stops, 30 communities, and 12 facilities) and 2206 arcs. The 18 transit lines each runs along one of the North/South or East/West roads. A total of 200 buses are in service with each route initially having between 6 and 14 vehicles and constant fleet size bounds of $y_l^{\min} = 1$ and $y_l^{\max} = 20$. Transit stops are located at each road intersection with some interspersed between intersections, and are used as OD pairs for day-to-day transit. Communities are distributed somewhat uniformly throughout the city, each having a population between 1000 and 3000. Facilities are distributed unevenly throughout the city with a constant quality metric of $S_j \equiv 1$.

Trials were set up to test various values of the model parameters $\beta$, $\mathcal{K}$, and $\epsilon$. For $\beta$ we tested $\beta \in \{0.5, 1.0, 1.5, 2.0\}$ while holding $\mathcal{K} = 6$ and $\epsilon = 0.01$. For $\mathcal{K}$ we tested $\mathcal{K} \in \{1, 6, 12, 18, 24, 30\}$ while holding $\beta = 1.0$ and $\epsilon = 0.01$. For $\epsilon$ we tested $\epsilon \in \{0.00, 0.01, 0.05, 0.10, \infty\}$ (with $\epsilon = \infty$ meaning that the user cost constraints were completely ignored) while holding $\beta = 1.0$ and $\mathcal{K} = 6$. For each parameter combination we solved the SAMP using 500 solution algorithm iterations followed by a neighborhood search of the best known solution. In order to reduce potential noise in the solutions from the random aspects of the search algorithm each problem instance was solved 10 times and the solution with the best objective value was kept. The results of these trials are discussed in Section \ref{subsec:smallscaleresults}.

% Results
\section{Computational Results}
\label{sec:results}

In this section we summarize the results of the computational trials discussed above. Complete data tables for all computational trials can be found in Appendix \ref{app:resultdata}.

%==============================================================================
\subsection{CTA Network Design Results}
\label{subsec:ctaresults}

\begin{figure}[h]
	\centering
	\subcaptionbox{Community area primary care gravity metrics under current CTA network.\label{subfig:metricsinitial}}{$\vcenter{\hbox{\includegraphics[height=0.25\textheight]{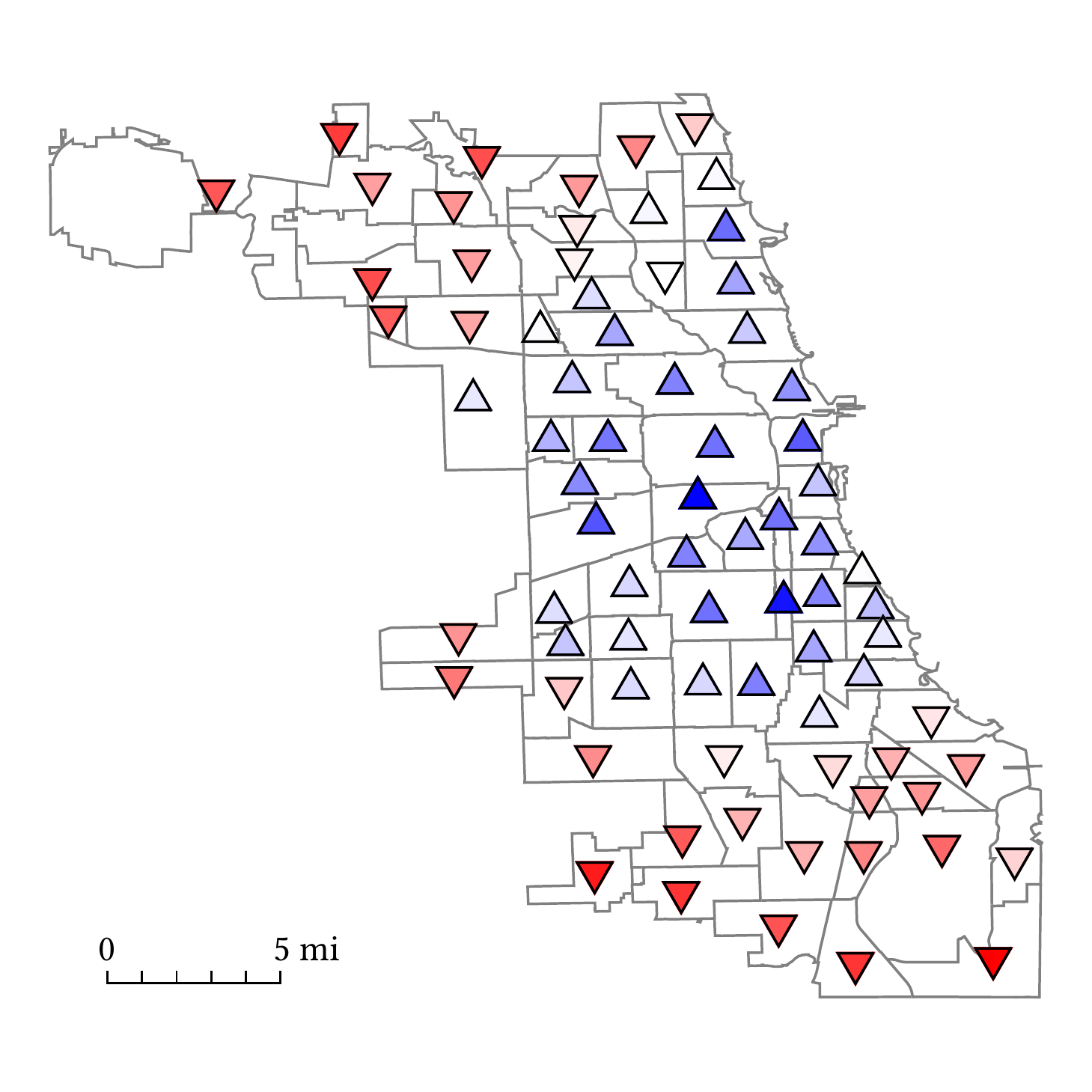}}}$} \qquad
	\subcaptionbox{Community area primary care gravity metrics after 500 search iterations.\label{subfig:metrics500}}{$\vcenter{\hbox{\includegraphics[height=0.25\textheight]{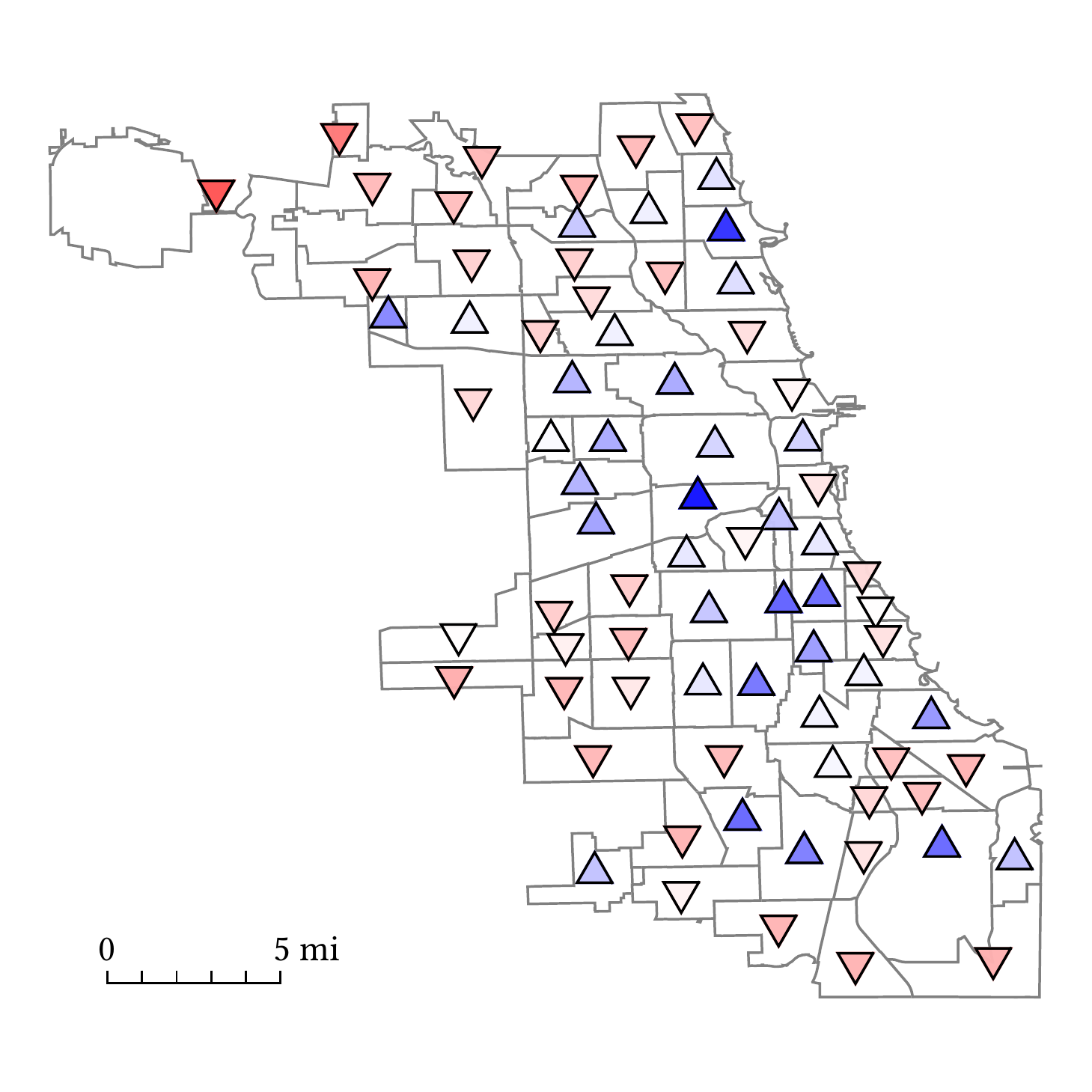}}}$} \qquad
	\begin{tabular}{c}
		\vspace{-2.3in} \\
		$\vcenter{\hbox{\includegraphics[height=0.13125\textheight]{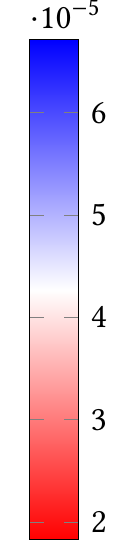}}}$
	\end{tabular}
	\caption{Maps comparing primary care access gravity metrics under the current CTA network and the solution after 500 search iterations. Each community area's metric is indicated by the color of the triangle located at its population-weighted centroid. The high and low ends of the color scale correspond to the highest and lowest initial gravity metrics in the city, respectively. The triangle points up if above the median and down if below.}
	\label{fig:metriccomparison}
\end{figure}

\begin{figure}[h]
	\centering
	$\vcenter{\hbox{\includegraphics[height=0.25\textheight]{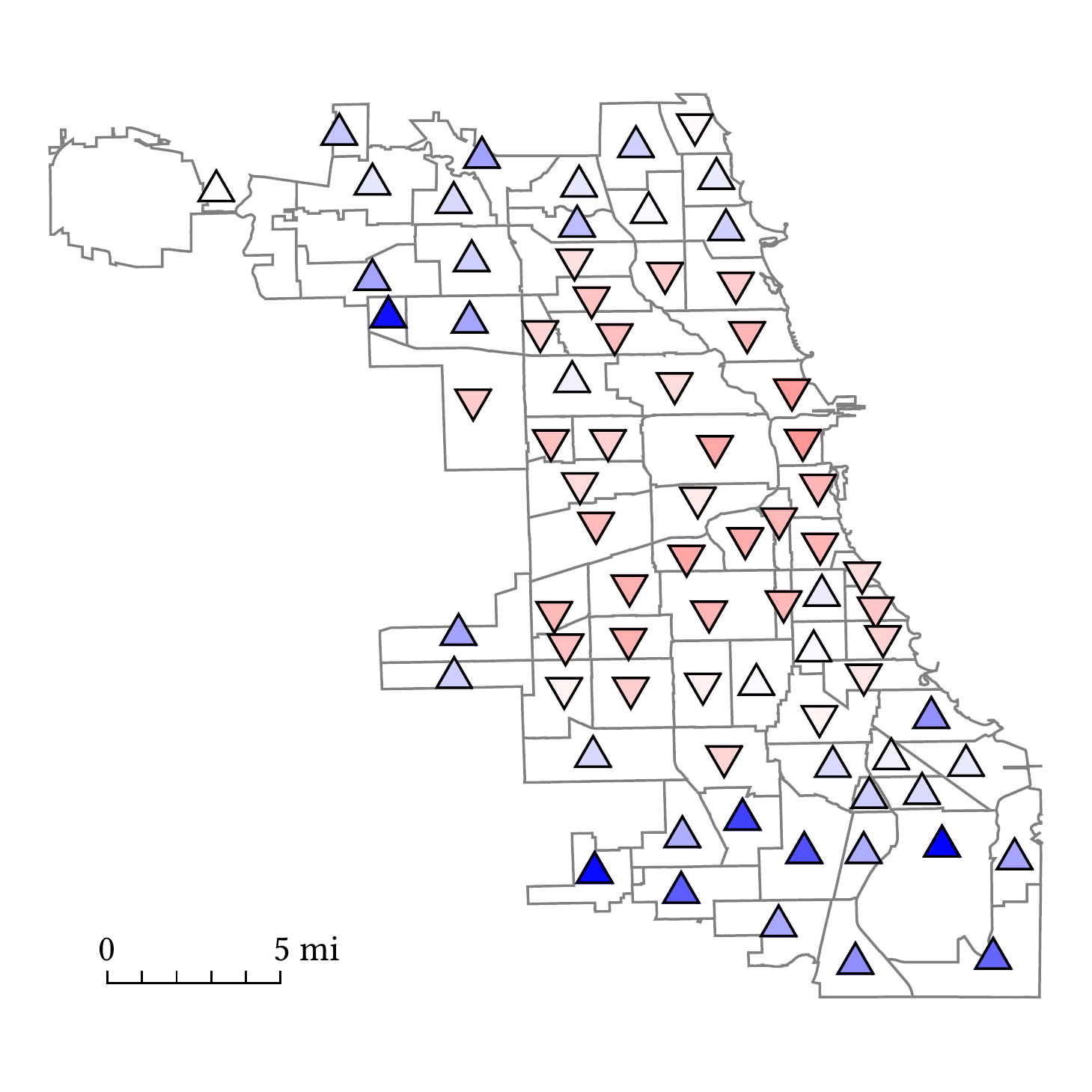}}}$ \qquad $\vcenter{\hbox{\includegraphics[height=0.13125\textheight]{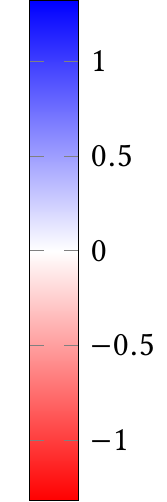}}}$
	\caption{Map of relative changes in primary care access gravity metrics after 500 search iterations. The exact center of the color scale is white and represents no change, with blue upward-pointing triangles indicating increase and red downward-pointing triangles indicating decrease. The high and low ends of the color scale correspond to $\pm 1.31949$, the magnitude of the greatest relative change of any community area (which was an increase for Mount Greenwood).}
	\label{fig:metricchange}
\end{figure}

The SAMP solution algorithm for the express design problem described in Section \ref{subsec:ctaproblem} was run for 500 iterations, after which the primary care accessibility metrics $A_i(\mathbf{y})$ were evaluated on the final network and compared to their original values. Table \ref{table:metricdatasummary} shows statistical results for the accessibility metrics of the 77 Chicago community areas. Of the 77 communities, approximately half (39 of 77) of the accessibility metrics improved while half (38 of 77) worsened. Among the 39 improved metrics the mean relative change was $+30.28551\%$, while among the 38 worsened metrics the mean relative change was $-12.27414\%$. The greatest individual increase among any community area's metric occurred for Mount Greenwood, which increased from $2.08703 \cdot 10^{-5}$ to $4.84086 \cdot 10^{-5}$ for a relative difference of $+131.94930\%$. The greatest individual decrease occurred for Near North Side, which decreased from $5.30391 \cdot 10^{-5}$ to $4.18952 \cdot 10^{-5}$ for a relative difference of $-21.01062\%$. Figure \ref{fig:metriccomparison} shows maps comparing the initial accessibility metric distribution to the final distribution, while Figure \ref{fig:metricchange} shows the relative difference in accessibility metrics before and after the trial. See Tables \ref{table:metricdata1}--\ref{table:metricdata2} in Appendix \ref{app:ctatables} for a complete list of metric changes over all community areas.

\begin{table}[h]
	\centering
	\begin{tabular}{l r r r r}
		\hline & \textbf{Initial} & \textbf{Final} & \textbf{Abs.\ Diff.} & \textbf{Rel.\ Diff.} \\
		\hline Mean & $4.18299 \cdot 10^{-5}$ & $4.31441 \cdot 10^{-5}$ & $+1.31417 \cdot 10^{-6}$ & $+3.14170\%$ \\
		Std.\ Dev.\ & $1.12720 \cdot 10^{-5}$ & $0.74786 \cdot 10^{-5}$ & $-3.79338 \cdot 10^{-6}$ & $-33.65312\%$ \\
		Median & $4.29677 \cdot 10^{-5}$ & $4.17583 \cdot 10^{-5}$ & $-1.20935 \cdot 10^{-6}$ & $-2.81456\%$ \\
		Max & $6.71548 \cdot 10^{-5}$ & $6.47152 \cdot 10^{-5}$ & $-2.43965 \cdot 10^{-6}$ & $-3.63287\%$ \\
		Min & $1.82806 \cdot 10^{-5}$ & $2.65753 \cdot 10^{-5}$ & $+8.29461 \cdot 10^{-6}$ & $+45.37372\%$ \\
		\hline
	\end{tabular}
	\caption{Summary statistics for the 77 Chicago community area accessibility metrics before and after running the express bus route design solution algorithm for 500 iterations.}
	\label{table:metricdatasummary}
\end{table}

\begin{figure}[h]
	\centering
	\includegraphics[width=0.45\textwidth]{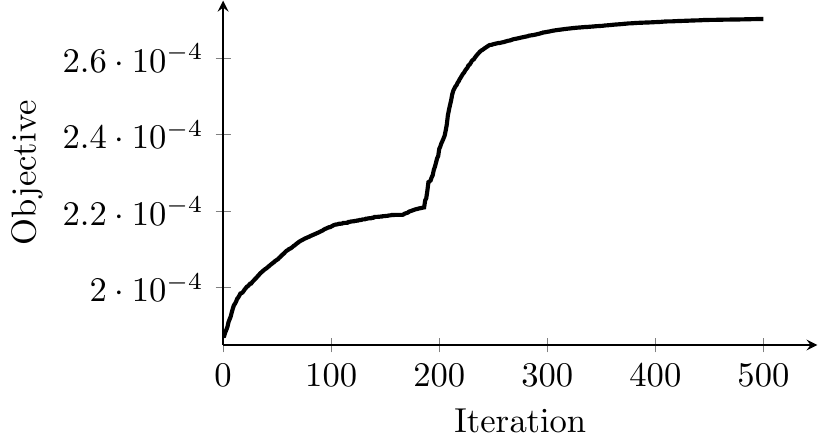}
	\caption{Objective value of the SAMP over 500 search iterations of the solution algorithm for the Chicago express route design problem.}
	\label{fig:tabusearchobjective}
\end{figure}

Figure \ref{fig:tabusearchobjective} shows the objective value of the express route design over the 500 search iterations. The objective increased from an initial value of $1.87302 \cdot 10^{-4}$ to a final value of $2.70226 \cdot 10^{-4}$, for an increase of $44.2\%$. The first 166 iterations consisted entirely of exchanging vehicles between non-express routes. The first express route was added during iteration 167, and the second was added during iteration 186, after which many more vehicles were added to express routes. This period of the search process corresponds to the sharp increase in the objective value seen near the middle of Figure \ref{fig:tabusearchobjective}. A total of 71 vehicles were diverted from existing lines to express lines ($4.25\%$ of the 1668 total buses in service), with 14 of the 65 available express routes receiving at least one vehicle.

Table \ref{table:expressfleets} in Appendix \ref{app:ctatables} shows the total number of vehicles added to each express line. Most received fewer than 3 vehicles, with the exceptions of the express runs of West 103rd (with 23 buses), Grand (with 22 buses), and Jeffery Local (with 8 buses). Tables \ref{table:positivefleets} and \ref{table:negativefleets} in Appendix \ref{app:ctatables} show the total numbers of vehicles added to or removed from existing lines, respectively. 13 lines received additional buses, with the most going to Vincennes/111th (gaining 72 buses), Harlem (gaining 71 buses), 31st (gaining 54 buses), South Michigan - OWL (gaining 27 buses), and Montrose (gaining 17 buses). 49 lines lost buses, 18 of which lost at least 5. The biggest losses occurred for 79th (losing 36 buses), Cottage Grove - OWL (losing 29 buses), King Drive (losing 26 buses), Western - OWL (losing 25 buses), and Halsted (losing 21 buses).

%==============================================================================
\subsection{CTA Network Supplemental Results}
\label{subsec:ctasupplemental}

The results of the supplemental tests discussed at the end of Section \ref{subsec:ctaproblem} are summarized below.

%%%
\paragraph{Express Line Exchanges:} A more restrictive version of the main express line design problem was run for which vehicles could not be exchanged freely between lines, but rather only exchanged between a line and its express run. The solution algorithm was run for 200 search iterations without finding a single solution that improved on the initial objective value. It was found that every possible SWAP move from the initial solution worsened the objective value, making the initial solution (with no exchanges) locally optimal.

%%%
\paragraph{Unconstrained Trial:} A relaxed version of the main express line design problem was run in which the user cost constraint (\ref{eqn:nlp1user}) was completely ignored. Eliminating this constraint also eliminated the need to evaluate the transit assignment model, drastically reducing the computation time needed for each search iteration. This trial was initialized with the final solution of the main trial and allowed to run for an additional 1000 iterations beyond the initial trial's 500 iterations. Between iterations 500 and 1500 the objective increased from $2.70226 \cdot 10^{-4}$ to $2.74911 \cdot 10^{-4}$, for an increase of $1.7\%$. The mean absolute value in relative difference of community accessibility metric was less than $2.1\%$. Beginning approximately with iteration 800 nearly every move involved adding a vehicle to Harlem, which ended with 413 vehicles. Vicennes/111th and 31st also ended with significantly more vehicles than the other bus lines (113 and 67, respectively), although neither total changed after approximately iteration 900.

%%%
\paragraph{Excluding Express Lines:} In order to determine whether the sharp improvement to the main trial's objective that occurred during iteration 186 was the result of the express lines, a restricted version of the trial was run which excluded all express lines. It was initialized with the main trial's solution at iteration 165, before the first express route was added, and allowed to run until reaching 250 iterations. Between iterations 165 and 250 this alternate trial's objective increased from $2.19014 \cdot 10^{-4}$ to $2.63517 \cdot 10^{-4}$, for a $0.5\%$ increase. During the same period the original trial had increased to $2.20147 \cdot 10^{-4}$, for a $20.3\%$ increase.

%==============================================================================
\subsection{Sensitivity Results from the Artificial Network}
\label{subsec:smallscaleresults}

Here we discuss the results of the different parameter tests separately. Table \ref{table:smallscalesol} in Appendix \ref{app:smallscaletables} shows the full solution vectors for all artificial network trials.

%%%
\paragraph{User Cost Increase Parameter:} In the case of $\epsilon = 0.00$ the solution algorithm made no changes to the initial network due to all neighbors of the initial solution being infeasible. All other values $\epsilon \in \{0.01,0.05,0.10,\infty\}$ produced exactly the same accessibility metrics for all communities, and thus exactly the same objective value, which increased from an initial value of $4.77960 \cdot 10^{-4}$ to a final value of $7.71440 \cdot 10^{-4}$ for a $61.4\%$ increase. Table \ref{table:epsilonchange} in Appendix \ref{app:smallscaletables} shows the full list of final community metrics. The solution vectors, themselves, were mostly the same across all trials. 12 of the 18 lines received the same number of vehicles in all solutions, with a maximum difference of 3 vehicles among the lines that did differ.

%%%
\paragraph{Community Inclusion Parameter:} The number of communities $\mathcal{K}$ to include in the objective function directly affects its value even if the solution remains the same. For this reason we cannot directly compare objective values for separate trials, but we can still compare the accessibility metrics for individual communities. Table \ref{table:krelchange} in Appendix \ref{app:smallscaletables} shows the full set of community access relative changes across the $\mathcal{K}$ trials.

The case of $\mathcal{K} = 1$ corresponds to optimizing the single lowest community area metric. This trial made no changes to the initial network, meaning that no explored solution was able to improve the accessibility of the lowest-metric community. Most of the larger values of $\mathcal{K}$ showed similar trends in which communities improved their access and which worsened, with the exception of $\mathcal{K} = 30$, which includes all communities in the objective (and thus corresponds to maximizing the average accessibility of all communities). Whereas the other trials generally improved the areas with the worst initial access at the cost of areas with the best initial access, the $\mathcal{K} = 30$ trial reversed this pattern and improved primarily the areas with the best initial access.

The initial artificial network contains 15 communities with below-median access and 15 with above-median access. Table \ref{table:kmedianmetrics} shows for each trial the average relative change for the initially below-median communities and the initially above-median communities. As mentioned above, $\mathcal{K} = 1$ resulted in no change. $\mathcal{K}$ between 6 and 24 all produced similar trends, achieving a small increase for the communities below median at the cost of a decrease for the communities above median. $\mathcal{K} = 30$ instead produces a small decrease for the communities below the median and a significant increase for the communities above the median.

\begin{table}[h]
	\centering
	\begin{tabular}{l r r}
		\hline \textbf{Trial} & \multicolumn{2}{c}{\textbf{Accessibility Rel.\ Diff.}} \\
		 & Init.\ Below Med.\ & Init.\ Above Med.\ \\
		\hline $\mathcal{K}=1$ & $0.00000\%$ & $0.00000\%$ \\
		$\mathcal{K}=6$ & $+1.94602\%$ & $-13.62515\%$ \\
		$\mathcal{K}=12$ & $+2.06421\%$ & $-15.13642\%$ \\
		$\mathcal{K}=18$ & $+1.70508\%$ & $-11.47495\%$ \\
		$\mathcal{K}=24$ & $+1.49110\%$ & $-9.81080\%$ \\
		$\mathcal{K}=30$ & $-1.04769\%$ & $+8.09278\%$ \\
		\hline
	\end{tabular}
	\caption{Average relative change in community accessibility metrics after artificial network trials for each value of $\mathcal{K}$. Each value represents the average relative change to the 15 communities initially above or below the median access level.}
	\label{table:kmedianmetrics}
\end{table}

%%%
\paragraph{Gravitational Decay Parameter:} The gravitational decay parameter $\beta$ changes the evaluations of all community accessibility metrics, and so neither the objective nor the individual metrics can be directly compared between separate trials. Relative changes can still be compared, and are shown in Table \ref{table:betarelchange} in Appendix \ref{app:smallscaletables}. All trials produced a qualitatively similar distribution of relative accessibility changes for each community, although the magnitudes of the relative differences tended to increase as $\beta$ increased. All trials produced the greatest relative increase and the greatest relative decrease for the same two communities. The rankings of the relative increases remained similar between trials, and 19 of the 30 communities showed the same sign of accessibility metric difference (increase or decrease) in all trials.

The trials for $\beta = 1.0$, $\beta = 1.5$, $\beta = 2.0$ produced nearly identical solution vectors, with 13 of the 18 lines receiving the same number of vehicles. Between $\beta = 1.0$ and $\beta = 1.5$ only 4 lines differed, all by 4 or fewer vehicles, and between $\beta = 1.5$ and $\beta = 2.0$ only 3 lines differed, all by 3 or fewer vehicles. The trial for $\beta = 0.5$ produced a significantly different solution vector from the rest, with 7 of the 18 fleet sizes differing from those of the other trials by 6 or more vehicles, 3 of which differed by 10 or more.

% Discussion, Conclusions, and Acknowledgements
\section{Discussion and Conclusions}
\label{sec:discussion}

The results presented in Section \ref{sec:results} are mixed but generally positive with respect to evaluating the effectiveness of the SAMP in meeting our modeling goals. While the main trial for express route design for the CTA network was able to significantly improve the primary care access levels of the most underserved areas of the city, the artificial network trials indicate that the model's output can be sensitive to certain parameters.

The results shown in Section \ref{subsec:ctaresults} indicate that the SAMP shows promise as a model for achieving better social equity via minor alterations to public transit. The main trial was able to achieve improved access levels for the least advantaged areas of the city. While this did require the access levels of the most advantaged areas to fall, the changes described in Figure \ref{fig:metriccomparison} and Table \ref{table:metricdatasummary} indicate that the relative gains far outweighed the relative losses. The overall distribution of accessibility metrics within the city remained relatively balanced, with approximately half improving and half worsening while the mean and median remained roughly the same, both changing by less than $3.2\%$. However, the spread of accessibility levels also narrowed significantly, with the standard deviation decreasing by more than $33.6\%$, the maximum decreasing by more than $3.6\%$, and the minimum increasing by more than $45.3\%$. Achieving a smaller spread in differential accessibility levels, with a particular focus on improving the worst case, is exactly the goal in improving equity of access, and in that regard the SAMP succeeded in its stated goal. Moreover, although the solution algorithm was halted after 500 search iterations, the relaxed run discussed in Section \ref{subsec:ctasupplemental} shows that relatively little improvement would be achieved if the algorithm had continued.

The form of the objective function $\mathrm{Access}(\mathbf{y})$ defined in (\ref{eqn:accesscost}) seems to have performed well. The resulting network generally improved the accessibility of not only the $\mathcal{K} = 8$ lowest-access areas but the surrounding areas as well. The artificial network results shown in Section \ref{subsec:smallscaleproblem}, however, indicate that the algorithm's performance is dependent on the choice of $\mathcal{K}$. While most intermediate values of $\mathcal{K}$ seem to result in qualitatively similar solutions, a very large choice of $\mathcal{K}$ causes no particular focus to be given to the least advantaged areas and may result in their access worsening, while a very small choice of $\mathcal{K}$ may lead to over-focusing on communities that cannot be significantly helped by any changes to the public transit system while ignoring other communities that can. This effect was also shown within the CTA study in the community area O'Hare, which is located at the westernmost bound of the Far North Side. O'Hare has relatively little CTA bus route access, and so none of our design decisions $\mathbf{y}$ could have a significant effect on the accessibility metric of O'Hare. By the end of the main trial its metric had increased by less than $+0.4\%$ in spite of being one of the $8$ lowest metrics (and thus part of the objective) throughout most of the algorithm's run time. The relaxed trial that ran for 1000 iterations after the initial trial had begun to divert a very large number of vehicles to the Harlem line, one of the few public transit connections for O'Hare, but was able to increase the objective by only $+1.7\%$.

The express run design component introduced in Section \ref{subsec:ctasupplemental} produced some unexpected results, and while the introduction of express runs seems to help with the accessibility objective, more significant system-wide changes are required before their effects are seen. The sharp improvement in the objective seen near iteration 200 of the main search in Figure \ref{fig:tabusearchobjective} did not occur in the trial without express routes available, indicating that they were essential for that improvement. However the trial in which only swaps between a route and its express run was unable to attain any improvement over the initial network. These results indicate that the express lines are not immediately helpful for the initial network, but that they do become helpful after significant alterations have been made to the network. Moreover, the total number of vehicles on express routes was fairly small, at only $4.25\%$ of the total bus fleet, which supports the idea that it is possible to significantly improve access levels with very little investment in new transit infrastructure.

The gravity decay metric trials described in Section \ref{subsec:smallscaleresults} indicate that the model's performance can depend greatly on the choice of parameters, echoing the results of the $\mathcal{K}$ trials. For most values of $\beta$ a similar set of communities showed an accessibility metric increase, indicating that all versions of the trial resulted in focus being given to a similar set of communities, however the magnitudes of the relative changes and the solution vectors were significantly different, indicating that the model is sensitive to how accessibility metrics are measured.

The tests of varying user cost increase bounds $\epsilon$ discussed in Section \ref{subsec:smallscaleresults} illustrate an important point regarding constraint (\ref{eqn:usercostbound}) of the SAMP. The case of $\epsilon = 0.00$ resulted in all neighbors of the initial solution being infeasible, which is unsurprising since the initial fleet sizes were chosen to be locally optimal with respect to user costs. For computational reasons the choice of $\epsilon = 0.00$ is not recommended for use. Due to the minor computational errors inherent in solving the transit assignment model \cite{spiess1989}, solutions whose user costs are very close to the upper bound $\mathrm{UserCost}(\mathbf{x^{\boldsymbol{*}}},\omega^*)$ may falsely be calculated as being slightly above it, marking those solutions erroneously as infeasible. In practice $\epsilon$ should be chosen to be at least a very small positive value larger than some reasonable numerical tolerance for the user cost calculation.

For the larger values of $\epsilon$ all trials produced identical objectives and community metrics, although not quite identical solution vectors (see Table \ref{table:smallscalesol}). The accessibility metric definition (\ref{eqn:gravitymetric})--(\ref{eqn:gravitycompetition}) is based on shortest path distance $d_{ij}(\mathbf{y})$ in $G$, and solution alterations that do not affect any community's shortest path have no effect on any community's accessibility metric. This is an undesirable property in an objective function since it produces regions of zero gradient within which there is no clear local search direction. Moreover, the fact that all values of $\epsilon$ produced nearly identical results indicates that the user cost constraint was not a significant impediment beyond the initial steps. This does not seem to be a realistic outcome for a real network, where it can be assumed that the initial network is already reasonably well optimized with respect to user cost.

The modeling framework of the SAMP shows promise for use in the problem of achieving social equity through use of minor public transit infrastructure changes. In the Chicago case study it was able to significantly improve the primary care access levels of underserved communities with relatively few changes to the current network. This represents just one of many potential applications of this type of model, but there is reason to be optimistic about the potential for achieving important social goals with very little investment. Due to the incompleteness of our data set and the number of simplifying assumptions made, the results of this study should not be regarded as a serious set of recommendations for changes in the CTA network design. In order to actually make such planning decisions much more data would be required and a more sophisticated user assignment model would be recommended. It would also be valuable to build the design of express runs into the model, itself, rather than simply generating a small number of candidates using heuristics.

The model presented in this study was meant to be general in its form, and the choices of objective and constraint function should be made according to the needs of the planner. While our access objective defined in (\ref{eqn:accesscost}) was generally able to drive the types of improvements we had hoped to see for the community areas, the sensitivity analysis for $\mathcal{K}$ indicates that this specific form could be improved. It would be reasonable, for example, to explore an objective that includes all accessibility metrics but dynamically weighted towards the current lowest values, or one that includes a set of the lowest-access communities alongside their neighboring communities.

%========================================================================
\section*{Acknowledgments}

The authors would like to thank Dr.~Lili Du, Dr.~Zongzhi Li, Dr.~Robert Ellis, and Dr.~Michael Pelsmajer for their valuable input. Thank you also to Dr.~Jan-Dirk Schm\"ocker for making their CapCon \cite{kurauchi2003} source code available for our use, and to Sharon Dunn for the computational resources used for generating our numerical results.

%==============================================================================
\small

% Use the \bibliography command to generate the .bbl file, then comment it out and switch to the \input command.

% Update bibliography with BibTeX
%\bibliographystyle{plainurl}
%\bibliography{references}

\begin{thebibliography}{10}

\bibitem{alhroob2014}
A.~Alhroob, H.~Tarawneh, and M.~Ayob.
\newblock A hybrid simulated annealing with tabu list and allowable solutions
  memory to solve university course timetabling problem.
\newblock {\em European Journal of Scientific Research}, 126(1):84--94, 2014.

\bibitem{arikan2012}
M.~Ar{\i}kan and S.~Erol.
\newblock A hybrid simulated annealing-tabu search algorithm for the part
  selection and machine loading problems in flexible manufacturing systems.
\newblock {\em International Journal of Advanced Manufacturing Technology},
  59:669--679, 2012.
\newblock \href {https://doi.org/10.1007/s00170-011-3506-0}
  {\path{doi:10.1007/s00170-011-3506-0}}.

\bibitem{bacharach1970}
M.~Bacharach.
\newblock {\em Biproportional Matrices and Input-Output Change}.
\newblock Cambridge University Press, London, 1970.

\bibitem{bagloee2011}
S.~A. Bagloee and A.~Ceder.
\newblock Transit-network design methodology for actual-size road networks.
\newblock {\em Transportation Research Part B: Methodological}, 45:1787--1804,
  2011.
\newblock \href {https://doi.org/10.1016/j.trb.2011.07.005}
  {\path{doi:10.1016/j.trb.2011.07.005}}.

\bibitem{healthcareunderserved}
J.~H. Beirne, E.~Salem, and R.~Ferguson.
\newblock Serving {C}hicago's underserved: {R}egional health systems profiles.
\newblock Technical report, Chicago Department of Public Health, Chicago Health
  \& Health Systems Project, Chicago, IL, 2005.

\bibitem{bellomo1970}
S.~J. Bellomo, R.~B. Dial, and A.~M. Voorhees.
\newblock Factors, trends, and guidelines related to trip length.
\newblock {\em National Cooperative Highway Research Program Report}, 89:1--61,
  1970.

\bibitem{linopt}
D.~Bertsimas and J.~N. Tsitsiklis.
\newblock {\em Introduction to Linear Optimization}.
\newblock Athena Scientific, Belmont, MA, 1997.

\bibitem{bielli2002}
M.~Bielli, M.~Caramia, and P.~Carotenuto.
\newblock Genetic algorithms in bus network optimization.
\newblock {\em Transportation Research Part C: Emerging Technologies},
  10:19--34, 2002.
\newblock \href {https://doi.org/10.1016/S0968-090X(00)00048-6}
  {\path{doi:10.1016/S0968-090X(00)00048-6}}.

\bibitem{buba2018}
A.~T. Buba and L.~S. Lee.
\newblock A differential evolution for simultaneous transit network design and
  frequency setting problem.
\newblock {\em Expert Systems With Applications}, 106:277--289, 2018.
\newblock \href {https://doi.org/10.1016/j.eswa.2018.04.011}
  {\path{doi:10.1016/j.eswa.2018.04.011}}.

\bibitem{cepeda2006}
M.~Cepeda, R.~Cominetti, and M.~Florian.
\newblock A frequency-based assignment model for congested transit networks
  with strict capacity constraints: {C}haracterization and computation of
  equilibria.
\newblock {\em Transportation Research Part B: Methodological}, 40:437--459,
  2006.
\newblock \href {https://doi.org/10.1016/j.trb.2005.05.006}
  {\path{doi:10.1016/j.trb.2005.05.006}}.

\bibitem{chakroborty2002}
P.~Chakroborty and T.~Dwivedi.
\newblock Optimal route network design for transit systems using genetic
  algorithms.
\newblock {\em Engineering Optimization}, 34(1):83--100, 2002.
\newblock \href {https://doi.org/10.1080/03052150210909}
  {\path{doi:10.1080/03052150210909}}.

\bibitem{cheung2016}
L.~L. Cheung and A.~S. Shalaby.
\newblock System-optimal re-routing transit assignment heuristic: {A}
  theoretical framework and large-scale case study.
\newblock In {\em Transportation Research Board 95th Annual Meeting}, number
  16--2333, 2016.
\newblock \href {https://doi.org/10.1016/j.ijtst.2017.09.001}
  {\path{doi:10.1016/j.ijtst.2017.09.001}}.

\bibitem{primarycare}
{Chicago Department of Public Health}.
\newblock {Public Health Services-- Chicago Primary Care Community Health
  Centers}, data last updated April 22, 2014.
\newblock URL:
  \url{https://data.cityofchicago.org/Health-Human-Services/Public-Health-Services-Chicago-Primary-Care-Commun/cjg8-dbka}.

\bibitem{traveltracker}
{Chicago Metropolitan Agency for Planning}.
\newblock {Travel Tracker Survey}, 2007 -- 2008, {C}hicago regional household
  travel inventory final report, data last updated October 7, 2014.
\newblock URL:
  \url{https://datahub.cmap.illinois.gov/dataset/traveltracker0708/resource/22eb7767-b143-416f-bb0f-37e6755231fd}.

\bibitem{ctabusridership}
{Chicago Transit Authority}.
\newblock {CTA} -- ridership -- avg.\ weekday bus stop boardings in {O}ctober
  2012, data last updated December 20, 2012.
\newblock URL:
  \url{https://data.cityofchicago.org/Transportation/CTA-Ridership-Avg-Weekday-Bus-Stop-Boardings-in-Oc/mq3i-nnqe}.

\bibitem{ctatrainridership}
{Chicago Transit Authority}.
\newblock {CTA} -- ridership -- `{L}' station entries -- monthly day-type
  averages \& totals, data last updated July 8, 2019.
\newblock URL:
  \url{https://data.cityofchicago.org/Transportation/CTA-Ridership-L-Station-Entries-Monthly-Day-Type-A/t2rn-p8d7}.

\bibitem{gtfs}
{Chicago Transit Authority}.
\newblock {GTFS}/scheduled service data, data last updated June 24, 2019.
\newblock URL: \url{https://www.transitchicago.com/developers/gtfs/}.

\bibitem{chu2018}
J.~C. Chu.
\newblock Mixed-integer programming model and branch-and-price-and-cut
  algorithm for urban bus network design and timetabling.
\newblock {\em Transportation Research Part B: Methodological}, 108:188--216,
  2018.
\newblock \href {https://doi.org/10.1016/j.trb.2017.12.013}
  {\path{doi:10.1016/j.trb.2017.12.013}}.

\bibitem{cipriani2012}
E.~Cipriani, S.~Gori, and M.~Petrelli.
\newblock Transit network design: {A} procedure and an application to a large
  urban area.
\newblock {\em Transportation Research Part C: Emerging Technologies},
  20:3--14, 2012.
\newblock \href {https://doi.org/10.1016/j.trc.2010.09.003}
  {\path{doi:10.1016/j.trc.2010.09.003}}.

\bibitem{delamater2013}
P.~L. Delamater.
\newblock Spatial accessibility in suboptimally configured health care systems:
  {A} modified two-step floating catchment area ({M2SFCA}) metric.
\newblock {\em Health \& Place}, 24:30--43, 2013.
\newblock \href {https://doi.org/10.1016/j.healthplace.2013.07.012}
  {\path{doi:10.1016/j.healthplace.2013.07.012}}.

\bibitem{fan2006}
W.~Fan and R.~B. Machemehl.
\newblock Optimal transit route network design problem with variable transit
  demand: {G}enetic algorithm approach.
\newblock {\em Journal of Transportation Engineering}, 132:40--51, 2006.
\newblock \href {https://doi.org/10.1061/(ASCE)0733-947X(2006)132:1(40)}
  {\path{doi:10.1061/(ASCE)0733-947X(2006)132:1(40)}}.

\bibitem{farahani2013}
R.~Z. Farahani, E.~Miandoabchi, W.~Y. Szeto, and H.~Rashidi.
\newblock A review of urban transportation network design problems.
\newblock {\em European Journal of Operational Research}, 229:281--302, 2013.
\newblock \href {https://doi.org/10.1016/j.ejor.2013.01.001}
  {\path{doi:10.1016/j.ejor.2013.01.001}}.

\bibitem{farber2014}
S.~Farber, M.~Z. Morang, and M.~J. Widener.
\newblock Temporal variability in transit-based accessibility to supermarkets.
\newblock {\em Applied Geography}, 53:149--159, 2014.
\newblock \href {https://doi.org/10.1016/j.apgeog.2014.06.012}
  {\path{doi:10.1016/j.apgeog.2014.06.012}}.

\bibitem{fu2012}
Q.~Fu, R.~Liu, and S.~Hess.
\newblock A review on transit assignment modelling approaches to congested
  networks: {A} new perspective.
\newblock {\em Procedia -- Social and Behavioral Sciences}, 54:1145--1155,
  2012.
\newblock \href {https://doi.org/10.1016/j.sbspro.2012.09.829}
  {\path{doi:10.1016/j.sbspro.2012.09.829}}.

\bibitem{healthcarepuzzle}
J.~Getzenberg, S.~Cohen, J.~Herd, J.~Sayer, and K.~Vandebroek.
\newblock The {C}hicago health care access puzzle: {F}itting the pieces
  together.
\newblock Technical report, Chicago Department of Public Health, Office of
  Policy \& Planning, Chicago, IL, 2008.

\bibitem{glover2013}
F.~Glover and M.~Laguna.
\newblock Tabu search: {E}ffective strategies for hard problems in analytics
  and computational science.
\newblock In P.~M. Pardalos, D.-Z. Du, and R.~Graham, editors, {\em Handbook of
  Combinatorial Optimization}, pages 3261--3362. Springer, Verlag, New York,
  2nd edition, 2013.

\bibitem{gu2010}
W.~Gu, X.~Wang, and S.~E. McGregor.
\newblock Optimization of preventive health care facility locations.
\newblock {\em International Journal of Health Geographics}, 9(1):1--16, 2010.
\newblock \href {https://doi.org/10.1186/1476-072X-9-17}
  {\path{doi:10.1186/1476-072X-9-17}}.

\bibitem{guihaire2008}
V.~Guihaire and J.-K. Hao.
\newblock Transit network design and scheduling: {A} global review.
\newblock {\em Transportation Research Part A: Policy and Practice},
  42:1251--1273, 2008.
\newblock \href {https://doi.org/10.1016/j.tra.2008.03.011}
  {\path{doi:10.1016/j.tra.2008.03.011}}.

\bibitem{joseph1982}
A.~E. Joseph and P.~R. Bantock.
\newblock Measuring potential physical accessibility to general practitioners
  in rural areas: {A} method and case study.
\newblock {\em Social Science \& Medicine}, 16(1):85--90, 1982.
\newblock \href {https://doi.org/10.1016/0277-9536(82)90428-2}
  {\path{doi:10.1016/0277-9536(82)90428-2}}.

\bibitem{kurauchi2003}
F.~Kurauchi, M.~G.~H. Bell, and J.-D. Schm{\"o}cker.
\newblock Capacity constrained transit assignment with common lines.
\newblock {\em Journal of Mathematical Modelling and Algorithms}, 2:309--327,
  2003.
\newblock \href {https://doi.org/10.1023/B:JMMA.0000020426.22501.c1}
  {\path{doi:10.1023/B:JMMA.0000020426.22501.c1}}.

\bibitem{lenin2016}
K.~Lenin, B.~R. Reddy, and M.~Suryakalavathi.
\newblock Hybrid tabu search-simulated annealing method to solve optimal
  reactive power problem.
\newblock {\em Electrical Power and Energy Systems}, 82:87--91, 2016.
\newblock \href {https://doi.org/10.1016/j.ijepes.2016.03.007}
  {\path{doi:10.1016/j.ijepes.2016.03.007}}.

\bibitem{luo2009}
W.~Luo and Y.~Qi.
\newblock An enhanced two-step floating catchment area ({E2SFCA}) method for
  measuring spatial accessibility to primary care physicians.
\newblock {\em Health \& Place}, 15:1100--1107, 2009.
\newblock \href {https://doi.org/10.1016/j.healthplace.2009.06.002}
  {\path{doi:10.1016/j.healthplace.2009.06.002}}.

\bibitem{luo2003b}
W.~Luo and F.~Wang.
\newblock Measures of spatial accessibility to health care in a {GIS}
  environment: {S}ynthesis and a case study in the {C}hicago region.
\newblock {\em Environment and Planning B: Planning and Design}, 30:865--884,
  2003.
\newblock \href {https://doi.org/10.1068/b29120} {\path{doi:10.1068/b29120}}.

\bibitem{luo2003a}
W.~Luo and F.~Wang.
\newblock Spatial accessibility to primary care and physician shortage area
  designation: {A} case study in {I}llinois with {GIS} aproaches.
\newblock In R.~Skinner and O.~Khan, editors, {\em Geographic Information
  Systems and Health Applications}, pages 260--278, Hershey, PA, 2003. Idea
  Group Publishing.
\newblock \href {https://doi.org/10.4018/978-1-59140-042-4.ch015}
  {\path{doi:10.4018/978-1-59140-042-4.ch015}}.

\bibitem{marcotte2004}
P.~Marcotte, S.~Nguyen, and A.~Schoeb.
\newblock A strategic flow model of traffic assignment in static capacitated
  networks.
\newblock {\em Operations Research}, 52(2):191--212, 2004.
\newblock \href {https://doi.org/10.1287/opre.1030.0091}
  {\path{doi:10.1287/opre.1030.0091}}.

\bibitem{mcdermot2017}
D.~McDermot, B.~Igoe, and M.~Stahre.
\newblock Assessment of healthy food availability in {W}ashington
  state--questioning of food desert paradigm.
\newblock {\em Journal of Nutrition Education and Behavior}, 49(2):130--136,
  2017.
\newblock \href {https://doi.org/10.1016/j.jneb.2016.10.012}
  {\path{doi:10.1016/j.jneb.2016.10.012}}.

\bibitem{mcgrail2009}
M.~R. McGrail and J.~S. Humphreys.
\newblock Measuring spatial accessibility to primary care in rural areas:
  {I}mproving the effectiveness of the two-step floating catchment area method.
\newblock {\em Applied Geography}, 29:533--541, 2009.
\newblock \href {https://doi.org/10.1016/j.apgeog.2008.12.003}
  {\path{doi:10.1016/j.apgeog.2008.12.003}}.

\bibitem{mcgrail2014}
M.~R. McGrail and J.~S. Humphreys.
\newblock Measuring spatial accessibility to primary health care services:
  {U}tilising dynamic catchment sizes.
\newblock {\em Applied Geography}, 54:182--188, 2014.
\newblock \href {https://doi.org/10.1016/j.apgeog.2014.08.005}
  {\path{doi:10.1016/j.apgeog.2014.08.005}}.

\bibitem{neutens2015}
T.~Neutens.
\newblock Accessibility, equity and health care: {R}eview and research
  directions for transport geographers.
\newblock {\em Journal of Transport Geography}, 43:14--27, 2015.
\newblock \href {https://doi.org/10.1016/j.jtrangeo.2014.12.006}
  {\path{doi:10.1016/j.jtrangeo.2014.12.006}}.

\bibitem{pearce2006}
J.~Pearce, K.~Witten, and P.~Bartie.
\newblock Neighbourhoods and health: {A} {GIS} approach to measuring community
  resource accessibility.
\newblock {\em Journal of Epidemiology \& Community Health}, 60(5):389--395,
  2006.
\newblock \href {https://doi.org/10.1136/jech.2005.04328}
  {\path{doi:10.1136/jech.2005.04328}}.

\bibitem{radke2000}
J.~Radke and L.~Mu.
\newblock Spatial decompositions, modeling and mapping service regions to
  predict access to social programs.
\newblock {\em Geographic Information Sciences}, 6(2):105--112, 2000.
\newblock \href {https://doi.org/10.1080/10824000009480538}
  {\path{doi:10.1080/10824000009480538}}.

\bibitem{p2codemaster}
A.~Rumpf.
\newblock Public transit optimization with social access objectives.
\newblock GitHub repository, 2020.
\newblock URL: \url{https://github.com/adam-rumpf/social-transit}.

\bibitem{p2data}
A.~Rumpf and H.~Kaul.
\newblock A public transit network optimization model for equitable access to
  social services.
\newblock Mendeley Data, 2021.
\newblock \href {https://doi.org/10.17632/pv2jvghs9b}
  {\path{doi:10.17632/pv2jvghs9b}}.

\bibitem{healthcaresafetynet}
E.~Salem and R.~Ferguson.
\newblock Casting {C}hicago's safety net: {A} 12-year review of {C}hicago's
  community-based primary care system.
\newblock Technical report, Chicago Department of Public Health, Planning and
  Development Division, Chicago, IL, 2005.

\bibitem{schmocker2011}
J.-D. Schm{\"o}cker, A.~Fonzone, H.~Shimamoto, F.~Kurauchi, and M.~G.~H. Bell.
\newblock Frequency-based transit assignment considering seat capacities.
\newblock {\em Transportation Research Part B: Methodological}, 45:392--408,
  2011.
\newblock \href {https://doi.org/10.1016/j.trb.2010.07.002}
  {\path{doi:10.1016/j.trb.2010.07.002}}.

\bibitem{schmocker2013}
J.-D. Schm{\"o}cker, H.~Shimamoto, and F.~Kurauchi.
\newblock Generation and calibration of transit hyperpaths.
\newblock {\em Procedia -- Social and Behavioral Sciences}, 80:211--230, 2013.
\newblock \href {https://doi.org/10.1016/j.sbspro.2013.05.013}
  {\path{doi:10.1016/j.sbspro.2013.05.013}}.

\bibitem{shafahi2018}
A.~Shafahi, Z.~Wang, and A.~Haghani.
\newblock {SpeedRoute}: {F}ast, efficient solutions for school bus routing
  problems.
\newblock {\em Transportation Research Part B: Methodological}, 117:473--493,
  2018.
\newblock \href {https://doi.org/10.1016/j.trb.2018.09.004}
  {\path{doi:10.1016/j.trb.2018.09.004}}.

\bibitem{spiess1989}
H.~Spiess and M.~Florian.
\newblock Optimal strategies: {A} new assignment model for transit networks.
\newblock {\em Transportation Research Part B: Methodological}, 23(2):83--102,
  1989.
\newblock \href {https://doi.org/10.1016/0191-2615(89)90034-9}
  {\path{doi:10.1016/0191-2615(89)90034-9}}.

\bibitem{wang2005}
F.~Wang and W.~Luo.
\newblock Assessing spatial and nonspatial factors for healthcare access:
  {T}owards an integrated approach to defining health professional shortage
  areas.
\newblock {\em Health \& Place}, 11:131--146, 2005.
\newblock \href {https://doi.org/10.1016/j.healthplace.2004.02.003}
  {\path{doi:10.1016/j.healthplace.2004.02.003}}.

\bibitem{weibull1976}
J.~W. Weibull.
\newblock An axiomatic approach to the measurement of accessibility.
\newblock {\em Regional Science and Urban Economics}, 6(4):357--379, 1976.
\newblock \href {https://doi.org/10.1016/0166-0462(76)90031-4}
  {\path{doi:10.1016/0166-0462(76)90031-4}}.

\end{thebibliography}

% Import predefined bibliography

%==============================================================================

\appendix

% Solution algorithm pseudocode
\section{Solution Algorithm Pseudocode}
\label{app:pseudocode}

The pseudocode for the hybrid tabu search/simulated annealing algorithm discussed in Section \ref{sec:algorithm} is included below. The source code for our C++ implementation can be viewed online \cite{p2codemaster}. For simplicity, $\mathbf{x}$ will be used to represent the entirety of the flow vector (although the assignment model also includes a waiting time scalar $\omega$). The canonical vector $\mathbf{e}_l \in \{0,1\}^{|L|}$ represents a vector with a $1$ in coordinate $l$ and $0$ elsewhere. Sets $\widetilde{N}^+$, $\widetilde{N}^-$, and $\widetilde{N}^\times$ represent the ADD, DROP, and SWAP neighbors kept as candidates after the first pass of the neighborhood search, respectively, while $N^+$, $N^-$, and $N^\times$ are the second-pass candidates. Bounds $\widetilde{N}_{\max}$ and $N_{\max}$ are the required numbers of candidates to be gathered during each pass. $Q_\text{in}^{\max}$ and $Q_\text{out}^{\max}$ are the inner and outer nonimprovement counter bounds, respectively.

\quad

\begin{algorithmic}[1]
\small
	\Procedure{// INITIALIZATION}{}
		\State $\mathbf{y}^{(0)} \gets \mathbf{y^{\boldsymbol{*}}}$ \Comment{current real-world fleet sizes}
		\State $\mathbf{y}_\text{best} \gets \mathbf{y}^{(0)}$ \Comment{best known solution and objective}
		\State $\mathrm{BestAccess} \gets \mathrm{Access}(\mathbf{y}^{(0)})$ \Comment{best known objective}
		\State $\mathrm{STM} \gets \emptyset$ \Comment{tabu move list, equipped with tabu tenures}
		\State $\mathrm{LTM} \gets \emptyset$ \Comment{attractive solution set}
		\State $T \gets T_0$ \Comment{simulated annealing temperature}
		\State $t \gets t_0$ \Comment{tabu tenure}
		\State $Q_{\text{in}} \gets 0$ \Comment{inner nonimprovement counter}
		\State $Q_{\text{out}} \gets 0$ \Comment{outer nonimprovement counter}
	\EndProcedure

	\For{$k = 0,1,\dots,K$} \Comment{search for a predetermined number $K$ of local moves}
		%%%
		\Procedure{// NEIGHBORHOOD SEARCH}{}
			%%%
			\State $\widetilde{N}^+ \gets \emptyset$ \Comment{tentative ADD neighborhood}
			\Repeat \Comment{first pass (ADD)}
				\State $l \gets \text{random element of $L$ not yet considered during this loop}$
				\If{$\mathbf{y}^{(k)} + \mathbf{e}_l$ satisfies vehicle limit bounds (\ref{eqn:nlp1vehicle})--(\ref{eqn:nlp1integer})} \Comment{consider only design-feasible moves}
					\If{$\mathbf{e}_l \notin \mathrm{STM}$ or $\mathrm{Access}(\mathbf{y}^{(k)} + \mathbf{e}_l) > \mathrm{BestAccess}$}
						\State $\widetilde{N}^+ \gets \widetilde{N}^+ \cup \{l\}$ \Comment{keep non-tabu candidates (unless improved best)}
					\EndIf
				\EndIf
			\Until{$|\widetilde{N}^+| \ge \widetilde{N}_{\text{max}}$ or all lines in $L$ have been considered}
			\State $N^+ \gets \emptyset$ \Comment{final ADD neighborhood}
			\Repeat \Comment{second pass (ADD)}
				\State $l \gets \text{unprocessed element of $\widetilde{N}^+$ with maximum $\mathrm{Access}(\mathbf{y}^{(k)} + \mathbf{e}_l)$}$
				\State $\mathbf{x}^{(k)} \gets \mathrm{TransitAssignment}(\mathbf{y}^{(k)}+\mathbf{e}_l)$ \Comment{calculate user flows for given move}
				\If{$(\mathbf{y}^{(k)}+\mathbf{e}_l,\mathbf{x}^{(k)})$ satisfy operator and user cost constraints (\ref{eqn:nlp1operator})--(\ref{eqn:nlp1user})}
					\State $N^+ \gets N^+ \cup \{l\}$ \Comment{keep only feasible candidates}
				\EndIf
			\Until{$|N^+| \ge N_\text{max}$ or all lines in $\widetilde{N}^+$ have been processed}
			%%%
			\State $\widetilde{N}^- \gets \emptyset$ \Comment{tentative DROP neighborhood}
			\Repeat \Comment{first pass (DROP)}
				\State $l \gets \text{random element of $L$ not yet considered during this loop}$
				\If{$\mathbf{y}^{(k)} - \mathbf{e}_l$ satisfies vehicle limit bounds (\ref{eqn:nlp1vehicle})--(\ref{eqn:nlp1integer})} \Comment{consider only design-feasible moves}
					\If{$-\mathbf{e}_l \notin \mathrm{STM}$ or $\mathrm{Access}(\mathbf{y}^{(k)} - \mathbf{e}_l) > \mathrm{BestAccess}$}
						\State $\widetilde{N}^- \gets \widetilde{N}^- \cup \{l\}$ \Comment{keep non-tabu candidates (unless improved best)}
					\EndIf
				\EndIf
			\Until{$|\widetilde{N}^-| \ge \widetilde{N}_{\text{max}}$ or all lines in $L$ have been considered}
			\State $N^- \gets \emptyset$ \Comment{final DROP neighborhood}
			\Repeat \Comment{second pass (DROP)}
				\State $l \gets \text{unprocessed element of $\widetilde{N}^-$ with maximum $\mathrm{Access}(\mathbf{y}^{(k)} - \mathbf{e}_l)$}$
				\State $\mathbf{x}^{(k)} \gets \mathrm{TransitAssignment}(\mathbf{y}^{(k)}-\mathbf{e}_l)$ \Comment{calculate user flows for given move}
				\If{$(\mathbf{y}^{(k)}-\mathbf{e}_l,\mathbf{x}^{(k)})$ satisfy operator and user cost constraints (\ref{eqn:nlp1operator})--(\ref{eqn:nlp1user})}
					\State $N^- \gets N^- \cup \{l\}$ \Comment{keep only feasible candidates}
				\EndIf
			\Until{$|N^-| \ge N_\text{max}$ or all lines in $\widetilde{N}^-$ have been processed}
			%%%
			\If{$|N^+| = |N^-| = 0$} \Comment{handle the event of a failed search}
				\If{bounds $\widetilde{N}_\text{max}$ were used in the previous neighborhood search}
					\State restart neighborhood search while ignoring the bounds $\widetilde{N}_\text{max}$
				\Else
					\State $Q_{\text{out}} \gets Q_{\text{out}}+1$ \Comment{work towards diversification}
					\State delete from $\mathrm{STM}$ the tabu rule with the shortest remaining tenure
					\State restart neighborhood search
				\EndIf
			\EndIf
			%%%
			\State $\widetilde{N}^\times \gets \set{(m,l) \in N^- \times N^+}{z_m = z_l}$ \Comment{set of all compatible ADD/SWAP pairs}
			\State $N^\times \gets \emptyset$ \Comment{final SWAP neighborhood}
			\Repeat
				\State $(m,l) \gets \text{unprocessed element of $N^\times$ with maximum $\mathrm{Access}(\mathbf{y}^{(k)} - \mathbf{e}_m) + \mathrm{Access}(\mathbf{y}^{(k)} + \mathbf{e}_l)$}$
				\State $\mathbf{x}^{(k)} \gets \mathrm{TransitAssignment}(\mathbf{y}^{(k)}-\mathbf{e}_m+\mathbf{e}_l)$
				\If{$(\mathbf{y}^{(k)}-\mathbf{e}_m+\mathbf{e}_l,\mathbf{x}^{(k)})$ satisfy operator and user cost constraints (\ref{eqn:nlp1operator})--(\ref{eqn:nlp1user})}
					\State $N^\times \gets N^\times \cup \{(m,l)\}$ \Comment{keep only feasible candidates}
				\EndIf
			\Until{$|N^\times| \ge N_\text{max}$ or all moves in $\widetilde{N}^\times$ have been processed}
			\State $\mathbf{\widetilde{y}}^{(k,1)} \gets \text{element of $N^+ \cup N^- \cup N^\times$ with greatest access}$ \Comment{returning two best neighbors}
			\State $\mathbf{\widetilde{y}}^{(k,2)} \gets \text{element of $N^+ \cup N^- \cup N^\times$ with second-greatest access}$
		\EndProcedure
		
		%%%
		\Procedure{// NEIGHBOR PROCESSING}{}
			\If{$\mathrm{Access}(\mathbf{\widetilde{y}}^{(k,1)}) - \mathrm{Access}(\mathbf{y}^{(k)}) > 0$}
				\State $Q_\text{out} \gets 0$ \Comment{improvement iteration}
				\State $t \gets t_0$
				\State $\mathbf{y}^{(k+1)} \gets \mathbf{\widetilde{y}}^{(k,1)}$ \Comment{move to improving neighbor}
				\State add rules with tenure $t$ to $\mathrm{STM}$ to prevent undoing the ADD/DROP moves of $\mathbf{y}^{(k+1)}$
				\If{$\mathrm{Access}(\mathbf{\widetilde{y}}^{(k,1)}) > \mathrm{BestAccess}$}
					\State $\mathbf{y}_\text{best} \gets \mathbf{\widetilde{y}}^{(k,1)}$ \Comment{update best known solution}
					\State $\mathrm{BestAccess} \gets \mathrm{Access}(\mathbf{\widetilde{y}}^{(k,1)})$
				\EndIf
			\Else
				\State $Q_\text{out} \gets Q_\text{out} + 1$ \Comment{nonimprovement iteration}
				\State $Q_\text{in} \gets Q_\text{in} + 1$
				\State $r \sim U[0,1]$ \Comment{random variable drawn uniformly from $[0,1]$}
				\If{$r < \exp\{-[ \mathrm{Access}(\mathbf{y}^{(k)}) - \mathrm{Access}(\mathbf{\widetilde{y}}^{(k,1)})]/T\}$}
					\State increment $t$ \Comment{SA criterion passed}
					\State $\mathbf{y}^{(k+1)} \gets \mathbf{\widetilde{y}}^{(k,1)}$ \Comment{move to non-improving neighbor}
					\State make undoing the ADD/DROP moves of $\mathbf{y}^{(k+1)}$ tabu with tenure $t$
					\State $Q_\text{in} \gets 0$
					\State $\mathrm{LTM} \gets \mathrm{LTM} \cup \{\mathbf{\widetilde{y}}^{(k,2)}\}$ \Comment{keep second-best neighbor as an attractive solution}
				\Else \Comment{SA criterion failed}
					\State $\mathrm{LTM} \gets \mathrm{LTM} \cup \{\mathbf{\widetilde{y}}^{(k,1)}\}$ \Comment{keep best neighbor as an attractive solution}
				\EndIf
			\EndIf
		\EndProcedure
		
		%%%
		\Procedure{// UPKEEP}{}
			\If{$Q_{\text{in}} \ge Q_{\text{in}}^{\text{max}}$}
				\State $Q_{\text{in}} \gets 0$ \Comment{attempt to diversify}
				\State $Q_{\text{out}} \gets Q_{\text{out}}+1$
				\State $\mathbf{y} \gets \text{randomly chosen element of $\mathrm{LTM}$}$
				\State $\mathrm{LTM} \gets \mathrm{LTM} \setminus \{\mathbf{y}\}$
				\State $\mathbf{y}^{(k+1)} \gets \mathbf{y}$ \Comment{move to a random attractive solution}
				\State increment $t$
			\EndIf
			\If{$Q_{\text{out}} \ge Q_{\text{out}}^{\text{max}}$}
				\State $t \gets t_0$ \Comment{attempt to intensify}
			\EndIf
			\State remove moves from $\mathrm{STM}$ whose tenures have elapsed
			\State apply cooling schedule to $T$
		\EndProcedure
	\EndFor

	%%%
	\State \textbf{return} $(\mathbf{y}_\text{best}, \mathrm{BestAccess})$
\end{algorithmic}

\quad

After the algorithm terminates and we obtain the best known solution, we conduct an exhaustive local search starting at $\mathbf{y}_\text{best}$. In each iteration of this local search we consider every possible ADD, DROP, and SWAP move without considering tabu rules, always taking the best neighbor and terminating when no neighbors offer any improvement. The final result is guaranteed to be locally optimal.

% CTA data processing
\section{CTA Network Data Processing}
\label{app:data}

This section gives a brief description of the data sources and processing used to generate the CTA computational trials described in Section \ref{subsec:ctaproblem}. The source code used for data processing can be viewed online \cite{p2codemaster}, as can the raw data used to define the network \cite{p2data}.

%==============================================================================
\subsection{Public Transit Network}
\label{app:networkdata}

Most of the structural information to construct transit network representation $G = (V,E)$ was obtained from the CTA General Transit Feed Specification (GTFS) files \cite{gtfs}. This includes a table of transit line stops (including geographic coordinates), transit lines, line trips, and the sequence of stops traversed on each line trip (including arrival/departure times). The final network included a total of 3656 nodes (including 997 stop and 2659 boarding) and 14,840 arcs (including 4912 line, 2659 boarding/alighting pairs, and 2305 walking).

%%%
\paragraph{Stop Clustering:} The GTFS files define over 11,000 stops in the CTA network. To reduce the size of the network a $k$-means process was used to partition the stops into clusters such that the total variance in position within all clusters was minimized. Each stop was then identified with the cluster with the nearest centroid, resulting in a set of 997 cluster centroids for use as stop nodes.  The mean geodesic distance (WGS 84 ellipsoid) between a stop and its associated cluster centroid was 0.13333 miles, with a standard deviation of 0.08327 miles, a median of 0.12103 miles, and a maximum of 0.84410 miles.

%%%
\paragraph{Transit Lines:} The GTFS stop time file was used to obtain the sequence of stops for each line $l \in L$, which were then mapped to the cluster centroids in the reduced network to define the stop nodes $V_\text{stop}^l$. Consecutive pairs of distinct stops defined line arcs $ij \in E_\text{line}^l$, with arrival time differences being used to define line arc travel times $c_{ij}$. If two consecutive stops had been identified with the same cluster centroid, then rather than generating a line arc as a loop, the travel time between the stops was equally distributed among the incident line arcs in order to preserve the total line travel time $\tau^l := \sum_{ij \in E_\text{line}^l} c_{ij}$. The total number of daily stop arrivals, divided by the number of unique stops, was used to calculate a daily average number of visits per stop $v_l$. The initial average line frequency was then calculated as $f_l^* := v_l/\tau^l$. This was used to estimate the initial fleet size as $y_l^* := \left\lceil \tau^l f_l^* \right\rceil$, obtained from the proportionality $f_l = y_l/\tau$ (see Section \ref{subsec:network}), taking the ceiling to ensure integrality.

The operational hours for the line were used to determine the fraction $\tau_f^l \in [0,1]$ of the full 1440 daily time horizon $\tau$ during which the line was active. During the execution of the solution algorithm the current value of $y_l$ was used to update the frequency and capacity of each line. Line capacities were scaled by a factor of $\tau_f^l$ and calculated as $u_l := \tau_f^l \tau b_l f_l$, the maximum number of people carried with service frequency $f_l$ over the actual operational hours of that line. Average line frequency was defined as $f_l := y_l/\tau^l$. The fleet size lower bound $y_l^{\min}$ for each bus line was taken as the minimum $y_l$ required to achieve a service frequency of at least one arrival per 30 minutes.

%%%
\paragraph{Walking Arcs:} Walking arcs $E_\text{walk}$ were generated between pairs of clustered stops less than $0.75$ miles from each other, measured using $L^1$-distance (a.k.a.\ taxicab distance or Manhattan distance) under the assumption that walking paths within city blocks are generally restricted to the cardinal directions. To prevent excessive numbers of arcs from being generated by dense cliques of pairwise-adjacent stops, walking arcs were not generated between a pair of stops if a path through an intermediate stop existed. For each pair of connected stops, a bidirectional pair of walking arcs was generated with a travel time calculated from the $L^1$-distance between the coordinates and an assumed walking speed of 4 ft/sec.

%==============================================================================
\subsection{Origin/Destination Travel Demand}
\label{app:oddemand}

An OD travel demand matrix is required by the public transit assignment model discussed in Section \ref{subsec:network}, and contains the travel demand between every pair of locations within the city via public transit during the day-to-day time horizon, estimated from ridership data. Because in general only stop-level data are available, the transit stops $V_\text{stop}$, themselves, are treated as the origins and destinations $V^+$ and $V^-$. The only available CTA stop-level bus data was from October 2012 \cite{ctabusridership}, and so stop-level train data from the same month were used \cite{ctatrainridership}. Each data set included total numbers of daily boardings. This lack of time-dependent data was the main reason behind using the entire 24-hour period as the time horizon $\tau$. Boardings at each stop were added to obtain an average total number of weekday boardings at each clustered stop. Alighting numbers were set equal to the boardings at each stop under the assumption that most daily trips are two-way.

These boarding and alighting numbers were used to estimate the OD matrix through iterative proportional fitting (IPF) \cite{bacharach1970}. The initial seed matrix was generated using a trip distribution function which assumed a gamma distribution of trip lengths \cite{bellomo1970}. Mean trip length was estimated using results from the 2007--2008 Chicago Metropolitan Agency for Planning (CMAP) Travel Tracker Survey \cite{traveltracker} which provides mean trip lengths for bus and train trips. These means were weighted by the total numbers of bus and train boardings to obtain a mode-weighted mean of $43.8$ minutes. The standard deviation was taken as $20$ minutes. The seed matrix was then generated by applying the resulting gamma distribution to the pairwise trip lengths in $G$. The IPF process stabilized after 10 iterations with an error of less than $0.001$, after which travel demands were rounded to the nearest integer.

%==============================================================================
\subsection{Community Nodes}
\label{app:communitydata}

Community nodes $\widetilde{V}^+$ represent population centers within Chicago, which were identified using 2010 census data. To reduce the size of the network, the over 800 census tracts were clustered into population-weighted centroids for each of the 77 Chicago community areas, grouping census tracts according to the community area to which they belonged. A community node was generated for each community area and connected via walking arcs to all stops within an $L^1$-distance of 1 mile. Exactly one population center (O'Hare) was located more than 1 mile from the nearest transit stop, and so its walking arc cutoff distance was extended to 2 miles. The total population over all census tracts was 2,739,860.

%==============================================================================
\subsection{Facility Nodes}
\label{app:facilitydata}

Facility nodes $\widetilde{V}^-$ correspond to primary care facilities within Chicago (see Figure \ref{fig:facility}). A list of primary care community health clinics and their geographic coordinates was obtained from the Chicago Data Portal \cite{primarycare}. These facilities all provide primary care, are open to everyone regardless of ability to pay, and have a specific mission to provide care for underserved populations. A facility node was generated for each primary care facility and connected via walking arcs to all stops within an $L^1$-distance of 1 mile.

%==============================================================================
\subsection{Candidate Express Route Generation}
\label{app:expressdata}

Candidate express routes for use in the design problem discussed in Section \ref{subsec:ctaproblem} were generated by applying a heuristic algorithm to all non-express CTA bus lines with at least 18 stops. Under the assumption that a beneficial express run would be one which makes it easier to reach a high-access area from a low-access area, a gravity metric $A_i(\mathbf{y^{\boldsymbol{*}}})$ as defined in (\ref{eqn:gravitymetric})--(\ref{eqn:gravitycompetition}) was calculated for all transit stops $i \in V_\text{stop}$ on the line in question, treating each stop as a community location with population $P_i \equiv 1$. The express run of the line then dropped all stops except for the 10\% highest-access and the 10\% lowest-access stops (rounded to the nearest integer). The line arcs of the express run had their travel times reduced by 40 seconds for each skipped stop. Express runs were assumed to have the same operation period $\tau_f^l$ as the original line. A total of 65 of the original 126 CTA bus lines were given an express run. These lines were used to define a modified version of network $G$ for use in some of the computational trials.

% Artificial network generation
\section{Artificial Network Generation}
\label{app:artificial}

In this section we describe the procedural network generation process used to generate the small artificial test case described in Section \ref{subsec:smallscaleproblem}. The network generation code and its resulting data can be viewed online \cite{p2codemaster,p2data}. The artificial network operated over a time horizon of $\tau = 60$ minutes, and was based on a $7 \times 11$ square lattice of roads representing a $0.9 \times 1.5$ mile city arranged in blocks. Transit stops were located on the intersections of the lattice with either 0, 1, or 2 additional stops (chosen uniformly at random) interspersed along the road segments between intersections. A total of 18 transit lines were defined, with one for each North/South and each East/West road on the grid. Arcs were assigned travel times based on Euclidean length and a bus speed chosen uniformly at random from $[0.3,0.6]$ miles per minute, plus a stopping time chosen uniformly at random from $[0.5,0.75]$ minutes. Trips began at one terminal, moved to the opposite terminal, and returned to the first terminal, followed by a layover chosen uniformly at random from $[20,30]$ minutes. Walking arcs connected each stop node to a maximum of 5 stop nodes falling within a cutoff $L^1$-distance of 0.12 miles, with travel times assuming a 4 ft/sec walking speed.

A total of 30 community nodes $\widetilde{V}^+$ were generated, with locations chosen uniformly at random within the city's $0.9 \times 1.5$ mile boundaries, and with the restriction that no two communities could be within 0.15 miles of each other in order to achieve somewhat even spacing. A population $P_i$ uniformly chosen from $[1000,3000]$ was assigned to each community $i$. A total of 12 facility nodes $\widetilde{V}^-$ were generated by again choosing coordinates uniformly at random, but with no minimum pairwise distance threshold, resulting in some densely-covered and some sparsely-covered areas of the city. All facilities $j$ were given a quality of $S_j \equiv 1$. Walking arcs were generated to connect communities and facilities to stops within 0.25 miles. The final network included 527 nodes and 2206 arcs.

An OD matrix of stop-level travel demand was generated using techniques similar to those discussed in Appendix \ref{app:oddemand}. The trip distribution function was a gamma distribution with a trip length mean of $2.75$ minutes and standard deviation of $1.25$ minutes. Stop-level total boarding numbers were determined by dividing 30\% of each community's population equally among all stops for which they were the nearest community. The seed matrix was perturbed by adding a random value chosen uniformly from $[-2,2]$ to each element, and then IPF was used to generate the OD matrix.

The initial fleet vector $\mathbf{y^{\boldsymbol{*}}}$ was generated by assigning a set total of 200 vehicles proportionally to each line according to its total travel demand, calculated as total OD demand for which that line was part of the shortest path through the network. This assignment was then followed with a local search of reassignment moves to minimize the user cost function (\ref{eqn:usercost}) under the assumption that the initial network should be reasonably-well optimized for that goal. The resulting fleet vector assigned between 6 and 14 vehicles to each line.

% Data tables
\section{Computational Trial Data Tables}
\label{app:resultdata}

This section includes the raw data tables for the computational trials described in Section \ref{sec:trials} and summarized in Section \ref{sec:results}. The raw data compiled in these tables can be viewed online \cite{p2data}.

%==============================================================================
\subsection{CTA Network Trial Data Tables}
\label{app:ctatables}

Tables \ref{table:metricdata1}--\ref{table:metricdata2} show the changes in all 77 Chicago community area accessibility metrics before and after running the express route design model for 500 iterations, and include the initial metric, final metric, absolute difference in metric, and relative different in metric before and after solving the model. The contents of these tables correspond to the summary statistics in Table \ref{table:metricdatasummary} and to Figures \ref{fig:metriccomparison}--\ref{fig:metricchange}.

Tables \ref{table:expressfleets}--\ref{table:negativefleets} show the changes in vehicle fleet sizes before and after running the express route design model for 500 iterations. Table \ref{table:expressfleets} shows the heuristically generated express lines described in Section \ref{subsec:ctaproblem}, along with the name of the line whose route they follow and the final number of vehicles. Table \ref{table:positivefleets} shows the original routes which received additional vehicles while \ref{table:negativefleets} shows the original routes which lost vehicles. Only nonzero changes are included in tables.

%==============================================================================
\subsection{Artificial Network Trial Data Tables}
\label{app:smallscaletables}

Tables \ref{table:epsilonchange}--\ref{table:smallscalesol} show the results of the artificial network trials described in Section \ref{subsec:smallscaleresults}. All tables list communities in the same order, from least to greatest initial accessibility metric. Table \ref{table:epsilonchange} shows the community metrics (along with absolute and relative differences) resulting from varying the allowable increase in user cost $\epsilon$. The cases of $\epsilon = 0.01$, $0.05$, $0.10$, and $\infty$ all produced identical metrics and so are listed together, while $\epsilon = 0.00$ produced no change and is not listed. Table \ref{table:krelchange} shows the relative change in community access metrics resulting from varying the number $\mathcal{K}$ of communities included in the objective. Note that the trial for $\mathcal{K}=1$ made no changes to the network, resulting in no relative change in accessibility metrics. Table \ref{table:betarelchange} shows the relative community access changes resulting from varying gravitational decay parameters $\beta$. The trials for $\beta = 1.5$ and $\beta = 2.0$ made no changes to the network. Table \ref{table:smallscalesol} shows the solution vectors for all artificial network trials. The columns correspond to each line's fleet size. Lines E0--E7 correspond to the East/West lines in order from South to North, while N0--N9 correspond to the North/South lines in order from West to East.

%==============================================================================

\FloatBarrier

\begin{table}[h]
	\centering
	\small
	\begin{tabular}{r l r r r r}
		\hline \multicolumn{2}{c}{\textbf{Community Area}} & \multicolumn{4}{c}{\textbf{Accessibility Metric}} \\
		ID & Name & Initial ($10^{-5}$) & Final ($10^{-5}$) & Abs.\ Diff.\ ($10^{-5}$) & Rel.\ Diff.\ \\
		\hline 1 & Rogers Park & $3.76050$ & $3.67836$ & $-0.08213$ & $-2.18412\%$ \\
		2 & West Ridge & $3.11850$ & $3.62795$ & $+0.50945$ & $+16.33651\%$ \\
		3 & Uptown & $5.67747$ & $6.18880$ & $+0.51133$ & $+9.00632\%$ \\
		4 & Lincoln Square & $4.34702$ & $4.41482$ & $+0.06780$ & $+1.55973\%$ \\
		5 & North Center & $4.27046$ & $3.68561$ & $-0.58485$ & $-13.69533\%$ \\
		6 & Lake View & $5.14559$ & $4.59176$ & $-0.55383$ & $-10.76320\%$ \\
		7 & Lincoln Park & $4.78835$ & $3.97534$ & $-0.81301$ & $-16.97888\%$ \\
		8 & Near North Side & $5.30391$ & $4.18952$ & $-1.11438$ & $-21.01062\%$ \\
		9 & Edison Park & $2.39304$ & $3.01503$ & $+0.62200$ & $+25.99205\%$ \\
		10 & Norwood Park & $3.35064$ & $3.62626$ & $+0.27562$ & $+8.22590\%$ \\
		11 & Jefferson Park & $3.24650$ & $3.67181$ & $+0.42531$ & $+13.10047\%$ \\
		12 & Forest Glen & $2.57881$ & $3.61349$ & $+1.03468$ & $+40.12241\%$ \\
		13 & North Park & $3.29882$ & $3.56305$ & $+0.26423$ & $+8.00975\%$ \\
		14 & Albany Park & $4.06956$ & $4.80290$ & $+0.73334$ & $+18.02014\%$ \\
		15 & Portage Park & $3.34545$ & $3.85255$ & $+0.50710$ & $+15.15776\%$ \\
		16 & Irving Park & $4.16694$ & $3.82487$ & $-0.34207$ & $-8.20907\%$ \\
		17 & Dunning & $2.56970$ & $3.57638$ & $+1.00668$ & $+39.17508\%$ \\
		18 & Montclare & $2.71115$ & $5.39024$ & $+2.67909$ & $+98.81737\%$ \\
		19 & Belmont Cragin & $3.41862$ & $4.41026$ & $+0.99164$ & $+29.00693\%$ \\
		20 & Hermosa & $4.29677$ & $3.82809$ & $-0.46868$ & $-10.90764\%$ \\
		21 & Avondale & $4.60055$ & $3.93870$ & $-0.66186$ & $-14.38643\%$ \\
		22 & Logan Square & $5.11390$ & $4.40383$ & $-0.71007$ & $-13.88517\%$ \\
		23 & Humboldt Park & $4.82359$ & $4.97242$ & $+0.14883$ & $+3.08539\%$ \\
		24 & West Town & $5.46999$ & $5.06966$ & $-0.40034$ & $-7.31875\%$ \\
		25 & Austin & $4.48806$ & $3.90469$ & $-0.58337$ & $-12.99818\%$ \\
		26 & West Garfield Park & $5.01456$ & $4.31704$ & $-0.69752$ & $-13.90987\%$ \\
		27 & East Garfield Park & $5.58700$ & $5.05871$ & $-0.52829$ & $-9.45562\%$ \\
		28 & Near West Side & $5.62957$ & $4.67567$ & $-0.95390$ & $-16.94447\%$ \\
		29 & North Lawndale & $5.39374$ & $4.98745$ & $-0.40629$ & $-7.53258\%$ \\
		30 & South Lawndale & $5.91802$ & $5.14372$ & $-0.77430$ & $-13.08373\%$ \\
		31 & Lower West Side & $6.71548$ & $6.47152$ & $-0.24396$ & $-3.63287\%$ \\
		32 & The Loop & $5.85658$ & $4.69402$ & $-1.16257$ & $-19.85058\%$ \\
		33 & Near South Side & $4.84068$ & $4.03246$ & $-0.80822$ & $-16.69641\%$ \\
		34 & Armour Square & $5.62103$ & $4.83526$ & $-0.78577$ & $-13.97914\%$ \\
		35 & Douglas & $5.31330$ & $4.49348$ & $-0.81981$ & $-15.42943\%$ \\
		36 & Oakland & $4.28715$ & $3.92353$ & $-0.36362$ & $-8.48164\%$ \\
		37 & Fuller Park & $6.53301$ & $5.72950$ & $-0.80351$ & $-12.29917\%$ \\
		38 & Grand Boulevard & $5.43644$ & $5.63635$ & $+0.19991$ & $+3.67722\%$ \\
		\hline
	\end{tabular}
	\caption{Accessibility metrics of Chicago community areas 1--38 before and after running the express bus route design solution algorithm for 500 iterations.}
	\label{table:metricdata1}
\end{table}

\begin{table}[p]
	\centering
	\small
	\begin{tabular}{r l r r r r}
		\hline \multicolumn{2}{c}{\textbf{Community Area}} & \multicolumn{4}{c}{\textbf{Accessibility Metric}} \\
		ID & Name & Initial ($10^{-5}$) & Final ($10^{-5}$) & Abs.\ Diff.\ ($10^{-5}$) & Rel.\ Diff.\ \\
		\hline 39 & Kenwood & $4.89197$ & $4.26635$ & $-0.62562$ & $-12.78865\%$ \\
		40 & Washington Park & $5.12451$ & $5.18802$ & $+0.06350$ & $+1.23922\%$ \\
		41 & Hyde Park & $4.48332$ & $3.98970$ & $-0.49362$ & $-11.01010\%$ \\
		42 & Woodlawn & $4.66761$ & $4.38398$ & $-0.28363$ & $-6.07656\%$ \\
		43 & South Shore & $4.02695$ & $5.24741$ & $+1.22046$ & $+30.30727\%$ \\
		44 & Chatham & $3.93504$ & $4.34573$ & $+0.41069$ & $+10.43668\%$ \\
		45 & Avalon Park & $3.52573$ & $3.69037$ & $+0.16463$ & $+4.66944\%$ \\
		46 & South Chicago & $3.32192$ & $3.57757$ & $+0.25566$ & $+7.69606\%$ \\
		47 & Burnside & $3.34628$ & $3.90227$ & $+0.55599$ & $+16.61519\%$ \\
		48 & Calumet Heights & $3.25548$ & $3.66385$ & $+0.40836$ & $+12.54389\%$ \\
		49 & Roseland & $3.52535$ & $5.46726$ & $+1.94191$ & $+55.08407\%$ \\
		50 & Pullman & $3.10465$ & $4.01719$ & $+0.91254$ & $+29.39260\%$ \\
		51 & South Deering & $2.82462$ & $5.66636$ & $+2.84174$ & $+100.60620\%$ \\
		52 & East Side & $3.85229$ & $4.85076$ & $+0.99846$ & $+25.91864\%$ \\
		53 & West Pullman & $2.61301$ & $3.58184$ & $+0.96883$ & $+37.07706\%$ \\
		54 & Riverdale & $2.32709$ & $3.57620$ & $+1.24911$ & $+53.67687\%$ \\
		55 & Hegewisch & $1.82806$ & $3.54108$ & $+1.71302$ & $+93.70672\%$ \\
		56 & Garfield Ridge & $3.23351$ & $4.25005$ & $+1.01654$ & $+31.43771\%$ \\
		57 & Archer Heights & $4.58707$ & $3.80274$ & $-0.78433$ & $-17.09875\%$ \\
		58 & Brighton Park & $4.65970$ & $3.78923$ & $-0.87047$ & $-18.68091\%$ \\
		59 & McKinley Park & $5.47325$ & $4.49003$ & $-0.98322$ & $-17.96407\%$ \\
		60 & Bridgeport & $5.09985$ & $4.17583$ & $-0.92401$ & $-18.11846\%$ \\
		61 & New City & $5.62979$ & $4.80549$ & $-0.82430$ & $-14.64183\%$ \\
		62 & West Elsdon & $4.82566$ & $4.13941$ & $-0.68624$ & $-14.22074\%$ \\
		63 & Gage Park & $4.52211$ & $3.65044$ & $-0.87167$ & $-19.27574\%$ \\
		64 & Clearing & $2.98121$ & $3.51580$ & $+0.53458$ & $+17.93177\%$ \\
		65 & West Lawn & $3.73206$ & $3.59777$ & $-0.13429$ & $-3.59825\%$ \\
		66 & Chicago Lawn & $4.61669$ & $4.06572$ & $-0.55098$ & $-11.93441\%$ \\
		67 & West Englewood & $4.64974$ & $4.49646$ & $-0.15329$ & $-3.29667\%$ \\
		68 & Englewood & $5.49465$ & $5.54644$ & $+0.05179$ & $+0.94255\%$ \\
		69 & Greater Grand Crossing & $4.51601$ & $4.38897$ & $-0.12703$ & $-2.81299\%$ \\
		70 & Ashburn & $3.16735$ & $3.61026$ & $+0.44291$ & $+13.98359\%$ \\
		71 & Auburn Gresham & $4.11289$ & $3.64952$ & $-0.46337$ & $-11.26623\%$ \\
		72 & Beverly & $2.69120$ & $3.57752$ & $+0.88631$ & $+32.93378\%$ \\
		73 & Washington Heights & $3.55078$ & $5.67823$ & $+2.12744$ & $+59.91481\%$ \\
		74 & Mount Greenwood & $2.08703$ & $4.84086$ & $+2.75383$ & $+131.94930\%$ \\
		75 & Morgan Park & $2.33348$ & $4.15515$ & $+1.82167$ & $+78.06660\%$ \\
		76 & O'Hare & $2.64728$ & $2.65753$ & $+0.01024$ & $+0.38687\%$ \\
		77 & Edgewater & $4.32469$ & $4.55498$ & $+0.23028$ & $+5.32486\%$ \\
		\hline
	\end{tabular}
	\caption{Accessibility metrics of Chicago community areas 39--77 before and after running the express bus route design solution algorithm for 500 iterations.}
	\label{table:metricdata2}
\end{table}

\begin{table}[p]
	\centering
	\small
	\begin{tabular}{r l r}
		\hline \multicolumn{2}{c}{\textbf{Transit Line}} & \textbf{Vehicles} \\
		ID & Name & Final \\
		\hline 103 & West 103rd & 23 \\
		65 & Grand & 22 \\
		15 & Jeffery Local & 8 \\
		81 & Lawrence - OWL & 3 \\
		62 & Archer - OWL & 3 \\
		67 & 67th-69th-71st & 2 \\
		30 & South Chicago & 2 \\
		18 & 16th-18th & 2 \\
		85 & Central & 1 \\
		84 & Peterson & 1 \\
		73 & Armitage & 1 \\
		44 & Wallace-Racine & 1 \\
		20 & Madison - OWL & 1 \\
		2 & Hyde Park Express & 1 \\
		\hline
	\end{tabular}
	\caption{Buses assigned to each candidate express route after running the express bus route design solution algorithm for 500 iterations. Line names indicate the existing line, while the vehicle number indicates the number of vehicles assigned to its express run. Only nonzero changes are shown.}
	\label{table:expressfleets}
\end{table}

\begin{table}[p]
	\centering
	\small
	\begin{tabular}{r l r r r r r}
		\hline \multicolumn{2}{c}{\textbf{Transit Line}} & \multicolumn{4}{c}{\textbf{Vehicles}} \\
		ID & Name & Initial & Final & Abs.\ Inc.\ & Rel.\ Inc.\ \\
		\hline 112 & Vincennes/111th & 8 & 80 & $+72$ & $+900.00000\%$ \\
		88 & Harlem & 4 & 75 & $+71$ & $+1775.00000\%$ \\
		30 & 31st & 13 & 67 & $+54$ & $+415.38462\%$ \\
		34 & South Michigan - OWL & 10 & 37 & $+27$ & $+270.00000\%$ \\
		78 & Montrose & 19 & 36 & $+17$ & $+89.47368\%$ \\
		90 & Austin & 12 & 27 & $+15$ & $+125.00000\%$ \\
		71 & 71st/South Shore & 18 & 32 & $+14$ & $+77.77778\%$ \\
		52A & South Kedzie & 10 & 19 & $+9$ & $+90.00000\%$ \\
		62H & Archer/Harlem & 4 & 9 & $+5$ & $+125.00000\%$ \\
		51 & 51st & 7 & 9 & $+2$ & $+28.57143\%$ \\
		49B & North Western & 8 & 10 & $+2$ & $+25.00000\%$ \\
		77 & Belmont - OWL & 32 & 33 & $+1$ & $+3.12500\%$ \\
		111 & 111th/King Drive & 9 & 10 & $+1$ & $+11.11111\%$ \\
		\hline
	\end{tabular}
	\caption{Positive changes to the fleet sizes of each existing line after running the express bus route design solution algorithm for 500 iterations. Only nonzero increases are shown.}
	\label{table:positivefleets}
\end{table}

\begin{table}[p]
	\centering
	\small
	\begin{tabular}{r l r r r r r}
		\hline \multicolumn{2}{c}{\textbf{Transit Line}} & \multicolumn{4}{c}{\textbf{Vehicles}} \\
		ID & Name & Initial & Final & Abs.\ Dec.\ & Rel.\ Dec.\ \\
		\hline 79 & 79th & 43 & 7 & $-36$ & $-83.72093\%$ \\
		4 & Cottage Grove - OWL & 35 & 6 & $-29$ & $-82.85714\%$ \\
		3 & King Drive & 32 & 6 & $-26$ & $-81.25000\%$ \\
		49 & Western - OWL & 38 & 13 & $-25$ & $-65.78947\%$ \\
		8 & Halsted & 34 & 13 & $-21$ & $-61.76471\%$ \\
		6 & Jackson Park Express & 23 & 4 & $-19$ & $-82.60870\%$ \\
		66 & Chicago - OWL & 36 & 20 & $-16$ & $-44.44444\%$ \\
		147 & Inner Drive/Michigan Express & 27 & 12 & $-15$ & $-55.55556\%$ \\
		J14 & Jeffery Jump & 20 & 5 & $-15$ & $-75.00000\%$ \\
		87 & Higgins & 22 & 9 & $-13$ & $-59.09091\%$ \\
		9 & Ashland - OWL & 40 & 27 & $-13$ & $-32.50000\%$ \\
		63 & 63rd - OWL & 27 & 15 & $-12$ & $-44.44444\%$ \\
		72 & North & 27 & 17 & $-10$ & $-37.03704\%$ \\
		62 & Archer - OWL & 23 & 15 & $-8$ & $-34.78261\%$ \\
		28 & Stony Island & 15 & 8 & $-7$ & $-46.66667\%$ \\
		55 & Garfield - OWL & 20 & 13 & $-7$ & $-35.00000\%$ \\
		52 & Kedzie/California & 27 & 22 & $-5$ & $-18.51852\%$ \\
		65 & Grand & 21 & 16 & $-5$ & $-23.80952\%$ \\
		5 & South Shore Night Bus - OWL & 8 & 4 & $-4$ & $-50.00000\%$ \\
		50 & Damen & 22 & 18 & $-4$ & $-18.18182\%$ \\
		53 & Pulaski - OWL & 29 & 25 & $-4$ & $-13.79310\%$ \\
		54B & South Cicero & 10 & 6 & $-4$ & $-40.00000\%$ \\
		67 & 67th-69th-71st & 19 & 15 & $-4$ & $-21.05263\%$ \\
		74 & Fullerton & 21 & 17 & $-4$ & $-19.04762\%$ \\
		80 & Irving Park & 21 & 17 & $-4$ & $-19.04762\%$ \\
		82 & Kimball-Homan & 27 & 23 & $-4$ & $-14.81481\%$ \\
		85 & North Central & 16 & 12 & $-4$ & $-25.00000\%$ \\
		10 & Museum of Science \& Industry & 7 & 4 & $-3$ & $-42.85714\%$ \\
		12 & Roosevelt & 20 & 17 & $-3$ & $-15.00000\%$ \\
		155 & Addison & 10 & 7 & $-3$ & $-30.00000\%$ \\ 
		20 & Madison - OWL & 25 & 22 & $-3$ & $-12.00000\%$ \\
		60 & Blue Island/26th - OWL & 22 & 19 & $-3$ & $-13.63636\%$ \\
		7 & Harrison & 7 & 4 & $-3$ & $-42.85714\%$ \\
		8A & South Halsted & 11 & 8 & $-3$ & $-27.27273\%$ \\
		106 & East 103rd & 5 & 3 & $-2$ & $-40.00000\%$ \\
		111A & Pullman Shuttle & 6 & 4 & $-2$ & $-33.33333\%$ \\
		115 & Pullman/115th & 11 & 9 & $-2$ & $-18.18182\%$ \\
		35 & 31st/35th & 17 & 15 & $-2$ & $-11.76471\%$ \\
		44 & Wallace-Racine & 10 & 8 & $-2$ & $-20.00000\%$ \\
		47 & 47th & 18 & 16 & $-2$ & $-11.11111\%$ \\
		92 & California/Dodge & 13 & 11 & $-2$ & $-15.38462\%$ \\
		95 & Central & 14 & 12 & $-2$ & $-14.28571\%$ \\
		1 & Bronzeville/Union Station & 3 & 2 & $-1$ & $-33.33333\%$ \\
		15 & Jeffery Local & 15 & 14 & $-1$ & $-6.66667\%$ \\
		29 & State & 23 & 22 & $-1$ & $-4.34783\%$ \\
		39 & Pershing & 7 & 6 & $-1$ & $-14.28571\%$ \\
		53A & South Pulaski & 14 & 13 & $-1$ & $-7.14286\%$ \\
		55N & 55th/Narragansett & 3 & 2 & $-1$ & $-33.33333\%$ \\
		X9 & Ashland Express & 12 & 11 & $-1$ & $-8.33333\%$ \\
		\hline
	\end{tabular}
	\caption{Negative changes to the fleet sizes of each existing line after running the express bus route design solution algorithm for 500 iterations. Only nonzero decreases are shown.}
	\label{table:negativefleets}
\end{table}

\begin{table}[p]
	\centering
	\small
	\begin{tabular}{r r r r r}
		\hline \multicolumn{2}{l}{\textbf{Community}} & \multicolumn{3}{c}{\textbf{Accessibility Metric}} \\
		ID & Init.\ ($10^{-4}$) & Final ($10^{-4}$) & Abs.\ Diff.\ ($10^{-4}$) & Rel.\ Diff.\ \\
		\hline 1 & 1.07653 & 1.18673 & $+0.11020$ & $+10.23670\%$ \\
		2 & 1.07761 & 1.10246 & $+0.02485$ & $+2.30643\%$ \\
		3 & 1.28453 & 1.28535 & $+0.00082$ & $+0.06377\%$ \\
		4 & 1.34093 & 1.34037 & $-0.00056$ & $-0.04187\%$ \\
		5 & 1.34462 & 1.40086 & $+0.05624$ & $+4.18293\%$ \\
		6 & 1.36099 & 1.40739 & $+0.04640$ & $+3.40944\%$ \\
		7 & 1.37413 & 1.45289 & $+0.07876$ & $+5.73177\%$ \\
		8 & 1.40414 & 1.40601 & $+0.00187$ & $+0.13326\%$ \\
		9 & 1.48147 & 1.49523 & $+0.01376$ & $+0.92850\%$ \\
		10 & 1.54621 & 1.56862 & $+0.02241$ & $+1.44942\%$ \\
		11 & 1.55761 & 1.54156 & $-0.01605$ & $-1.03066\%$ \\
		12 & 1.60705 & 1.65760 & $+0.05055$ & $+3.14540\%$ \\
		13 & 1.63913 & 1.60793 & $-0.03120$ & $-1.90326\%$ \\
		14 & 1.64878 & 1.62623 & $-0.02255$ & $-1.36754\%$ \\
		15 & 1.66128 & 1.62850 & $-0.03278$ & $-1.97332\%$ \\
		16 & 1.66200 & 1.61250 & $-0.04951$ & $-2.97875\%$ \\
		17 & 1.66257 & 1.64956 & $-0.01301$ & $-0.78257\%$ \\
		18 & 1.71999 & 1.72415 & $+0.00416$ & $+0.24199\%$ \\
		19 & 1.72518 & 1.72665 & $+0.00146$ & $+0.08472\%$ \\
		20 & 1.73728 & 1.77214 & $+0.03486$ & $+2.00660\%$ \\
		21 & 1.75121 & 1.76912 & $+0.01792$ & $+1.02305\%$ \\
		22 & 1.76028 & 1.75492 & $-0.00536$ & $-0.30472\%$ \\
		23 & 1.78722 & 1.76766 & $-0.01956$ & $-1.09432\%$ \\
		24 & 2.04493 & 2.01483 & $-0.03010$ & $-1.47188\%$ \\
		25 & 2.10664 & 2.11388 & $+0.00724$ & $+0.34389\%$ \\
		26 & 2.22927 & 2.12694 & $-0.10233$ & $-4.59012\%$ \\
		27 & 2.30707 & 2.25028 & $-0.05679$ & $-2.46154\%$ \\
		28 & 2.70816 & 2.66493 & $-0.04323$ & $-1.59644\%$ \\
		29 & 3.15152 & 3.11043 & $-0.04109$ & $-1.30373\%$ \\
		30 & 3.32601 & 3.30135 & $-0.02466$ & $-0.74133\%$ \\
		\hline
	\end{tabular}
	\caption{Accessibility metrics for $\epsilon \in \{0.01,0.05,0.10,\infty\}$, all of which produced identical metrics. The case of $\epsilon = 0.00$ produced no change and is not shown.}
	\label{table:epsilonchange}
\end{table}

\begin{table}[p]
	\centering
	\small
	\begin{tabular}{r r r r r r r r}
		\hline \multicolumn{2}{l}{\textbf{Community}} & \multicolumn{6}{c}{\textbf{Accessibility Metric Relative Difference}} \\
		ID & Init.\ ($10^{-4}$) & $\mathcal{K}=1$ & $\mathcal{K}=6$ & $\mathcal{K}=12$ & $\mathcal{K}=18$ & $\mathcal{K}=24$ & $\mathcal{K}=30$ \\
		\hline 1 & 1.07653 & $0.00000\%$ & $+10.23670\%$ & $+11.52970\%$ & $+7.48895\%$ & $+7.75733\%$ & $-11.95180\%$ \\
		2 & 1.07761 & $0.00000\%$ & $+2.30643\%$ & $+2.21118\%$ & $+2.50718\%$ & $+2.26418\%$ & $+0.22101\%$ \\
		3 & 1.28453 & $0.00000\%$ & $+0.06377\%$ & $-0.17949\%$ & $+0.00831\%$ & $-0.17967\%$ & $-1.00109\%$ \\
		4 & 1.34093 & $0.00000\%$ & $-0.04187\%$ & $-0.17296\%$ & $-0.01994\%$ & $-0.25074\%$ & $-2.11496\%$ \\
		5 & 1.34462 & $0.00000\%$ & $+4.18293\%$ & $+4.12391\%$ & $+4.57557\%$ & $+4.34197\%$ & $+2.92032\%$ \\
		6 & 1.36099 & $0.00000\%$ & $+3.40944\%$ & $+6.03928\%$ & $+2.97023\%$ & $+2.23649\%$ & $-0.48357\%$ \\
		7 & 1.37413 & $0.00000\%$ & $+5.73177\%$ & $+5.50925\%$ & $+5.67805\%$ & $+5.48597\%$ & $-3.83201\%$ \\
		8 & 1.40414 & $0.00000\%$ & $+0.13326\%$ & $-0.04220\%$ & $+0.08314\%$ & $-0.09552\%$ & $-0.04932\%$ \\
		9 & 1.48147 & $0.00000\%$ & $+0.92850\%$ & $+0.60524\%$ & $+0.90231\%$ & $+0.75014\%$ & $-2.80432\%$ \\
		10 & 1.54621 & $0.00000\%$ & $+1.44942\%$ & $+1.24680\%$ & $+1.29418\%$ & $+1.03556\%$ & $-4.35472\%$ \\
		11 & 1.55761 & $0.00000\%$ & $-1.03066\%$ & $-1.23426\%$ & $-1.12582\%$ & $-1.33133\%$ & $+10.88780\%$ \\
		12 & 1.60705 & $0.00000\%$ & $+3.14540\%$ & $+2.86421\%$ & $+3.10455\%$ & $+2.94187\%$ & $-4.74928\%$ \\
		13 & 1.63913 & $0.00000\%$ & $-1.90326\%$ & $-2.10939\%$ & $-2.04352\%$ & $-2.26721\%$ & $-0.21303\%$ \\
		14 & 1.64878 & $0.00000\%$ & $-1.36754\%$ & $-1.49240\%$ & $-1.55215\%$ & $-1.81234\%$ & $+2.85737\%$ \\
		15 & 1.66128 & $0.00000\%$ & $-1.97332\%$ & $-2.09528\%$ & $-1.45453\%$ & $-0.71927\%$ & $+3.14082\%$ \\
		16 & 1.66200 & $0.00000\%$ & $-2.97875\%$ & $-3.12205\%$ & $+0.74136\%$ & $+1.98271\%$ & $-2.06459\%$ \\
		17 & 1.66257 & $0.00000\%$ & $-0.78257\%$ & $-1.27512\%$ & $-1.38717\%$ & $-1.69549\%$ & $+2.29955\%$ \\
		18 & 1.71999 & $0.00000\%$ & $+0.24199\%$ & $+0.12119\%$ & $+0.46545\%$ & $+0.57221\%$ & $-0.77614\%$ \\
		19 & 1.72518 & $0.00000\%$ & $+0.08472\%$ & $-0.03743\%$ & $+1.32173\%$ & $+2.26348\%$ & $-0.64065\%$ \\
		20 & 1.73728 & $0.00000\%$ & $+2.00660\%$ & $+1.82355\%$ & $+1.82582\%$ & $+1.53470\%$ & $-3.40562\%$ \\
		21 & 1.75121 & $0.00000\%$ & $+1.02305\%$ & $+0.85984\%$ & $+0.98778\%$ & $+0.80087\%$ & $-1.83097\%$ \\
		22 & 1.76028 & $0.00000\%$ & $-0.30472\%$ & $+0.54420\%$ & $-0.60213\%$ & $-0.80462\%$ & $-1.47586\%$ \\
		23 & 1.78722 & $0.00000\%$ & $-1.09432\%$ & $-0.47427\%$ & $-0.60620\%$ & $-0.73652\%$ & $+2.72516\%$ \\
		24 & 2.04493 & $0.00000\%$ & $-1.47188\%$ & $-1.95233\%$ & $-2.31021\%$ & $-1.34374\%$ & $+3.07357\%$ \\
		25 & 2.10664 & $0.00000\%$ & $+0.34389\%$ & $-0.49598\%$ & $-0.99791\%$ & $-1.31535\%$ & $+1.44756\%$ \\
		26 & 2.22927 & $0.00000\%$ & $-4.59012\%$ & $-4.67559\%$ & $-3.56198\%$ & $-3.26678\%$ & $+1.83404\%$ \\
		27 & 2.30707 & $0.00000\%$ & $-2.46154\%$ & $-2.54728\%$ & $-2.64045\%$ & $-2.63787\%$ & $+1.71891\%$ \\
		28 & 2.70816 & $0.00000\%$ & $-1.59644\%$ & $-1.78781\%$ & $-1.65567\%$ & $-1.77671\%$ & $+2.98058\%$ \\
		29 & 3.15152 & $0.00000\%$ & $-1.30373\%$ & $-1.37948\%$ & $-1.54499\%$ & $-1.73567\%$ & $+1.03780\%$ \\
		30 & 3.32601 & $0.00000\%$ & $-0.74133\%$ & $-0.73786\%$ & $-1.51039\%$ & $-1.65202\%$ & $+1.16945\%$ \\
		\hline
	\end{tabular}
	\caption{Relative changes in community metrics for different values of $\mathcal{K}$. Relative changes indicate difference between the final output of the solution algorithm and the initial metric for that community.}
	\label{table:krelchange}
\end{table}

\begin{table}[p]
	\centering
	\small
	\begin{tabular}{r r r r r}
		\hline \multicolumn{2}{l}{\textbf{Community}} & \multicolumn{3}{c}{\textbf{Accessibility Metric Rel.\ Diff.}} \\
		ID & $\beta = 0.5$ & $\beta = 1.0$ & $\beta = 1.5$ & $\beta = 2.0$ \\
		\hline 1 & $+5.02349\%$ & $+10.23670\%$ & $+17.66870\%$ & $+27.54780\%$ \\
		2 & $+1.15748\%$ & $+2.30643\%$ & $+3.29771\%$ & $+3.98499\%$ \\
		3 & $+0.12525\%$ & $+0.06377\%$ & $-0.45141\%$ & $-1.40403\%$ \\
		4 & $-0.20341\%$ & $-0.04187\%$ & $+0.13229\%$ & $-0.06179\%$ \\
		5 & $+2.63748\%$ & $+4.18293\%$ & $+4.64409\%$ & $+4.32943\%$ \\
		6 & $+2.51144\%$ & $+3.40944\%$ & $+0.13601\%$ & $-2.39743\%$ \\
		7 & $+2.20261\%$ & $+5.73177\%$ & $+10.51130\%$ & $+17.15080\%$ \\
		8 & $+0.09534\%$ & $+0.13326\%$ & $-0.91723\%$ & $-3.30634\%$ \\
		9 & $+0.78419\%$ & $+0.92850\%$ & $-0.18574\%$ & $-0.68032\%$ \\
		10 & $+0.03285\%$ & $+1.44942\%$ & $+3.48831\%$ & $+6.81105\%$ \\
		11 & $-0.98510\%$ & $-1.03066\%$ & $-1.44343\%$ & $-2.22022\%$ \\
		12 & $+1.53500\%$ & $+3.14540\%$ & $+4.46802\%$ & $+5.46441\%$ \\
		13 & $-1.29035\%$ & $-1.90326\%$ & $-1.83503\%$ & $-1.42341\%$ \\
		14 & $-0.97666\%$ & $-1.36754\%$ & $-1.14218\%$ & $-0.34302\%$ \\
		15 & $-1.32526\%$ & $-1.97332\%$ & $-1.62264\%$ & $-1.07859\%$ \\
		16 & $-2.05815\%$ & $-2.97875\%$ & $-2.76270\%$ & $-1.87015\%$ \\
		17 & $-0.24560\%$ & $-0.78257\%$ & $-0.35359\%$ & $+0.59357\%$ \\
		18 & $-0.52186\%$ & $+0.24199\%$ & $+2.54721\%$ & $+6.52475\%$ \\
		19 & $-0.48299\%$ & $+0.08472\%$ & $+1.03859\%$ & $+2.05805\%$ \\
		20 & $+0.39778\%$ & $+2.00660\%$ & $+3.79349\%$ & $+6.22916\%$ \\
		21 & $+0.49524\%$ & $+1.02305\%$ & $+1.01055\%$ & $-0.19708\%$ \\
		22 & $-0.26372\%$ & $-0.30472\%$ & $+0.19110\%$ & $+0.98391\%$ \\
		23 & $+0.12488\%$ & $-1.09432\%$ & $-1.27936\%$ & $-1.31398\%$ \\
		24 & $-0.83268\%$ & $-1.47188\%$ & $-1.30115\%$ & $-0.89414\%$ \\
		25 & $+0.44535\%$ & $+0.34389\%$ & $+0.25383\%$ & $+0.10043\%$ \\
		26 & $-2.24379\%$ & $-4.59012\%$ & $-5.91311\%$ & $-7.05280\%$ \\
		27 & $-1.56020\%$ & $-2.46154\%$ & $-2.23398\%$ & $-1.84788\%$ \\
		28 & $-1.07935\%$ & $-1.59644\%$ & $-1.40627\%$ & $-0.78129\%$ \\
		29 & $-0.91016\%$ & $-1.30373\%$ & $-1.28732\%$ & $-1.00115\%$ \\
		30 & $-0.35668\%$ & $-0.74133\%$ & $-0.37780\%$ & $-0.11457\%$ \\
		\hline
	\end{tabular}
	\caption{Relative changes in community metrics for different values of $\beta$. Relative changes indicate difference between the final output of the solution algorithm and the initial metric for that community, both measured using the corresponding value of $\beta$.}
	\label{table:betarelchange}
\end{table}

\begin{table}[p]
	\centering
	\footnotesize
	\begin{tabular}{r l r r r r r r r r r r r r r r r r r r}
		\hline \multicolumn{2}{c}{\textbf{Trial}} & \multicolumn{18}{c}{\textbf{Solution Vector}} \\
		 & & E0 & E1 & E2 & E3 & E4 & E5 & E6 & E7 & N0 & N1 & N2 & N3 & N4 & N5 & N6 & N7 & N8 & N9 \\
		\hline & Initial & 8 & 12 & 12 & 11 & 12 & 14 & 11 & 12 & 12 & 10 & 13 & 13 & 12 & 11 & 6 & 7 & 12 & 12 \\
		\hline & $\epsilon = 0.00$ & 8 & 12 & 12 & 11 & 12 & 14 & 11 & 12 & 12 & 10 & 13 & 13 & 12 & 11 & 6 & 7 & 12 & 12 \\
		$*$ & $\epsilon = 0.01$ & 1 & 16 & 4 & 8 & 12 & 14 & 11 & 12 & 12 & 10 & 20 & 10 & 7 & 20 & 6 & 5 & 20 & 12 \\
		 & $\epsilon = 0.05$ & 1 & 16 & 2 & 11 & 12 & 14 & 11 & 12 & 12 & 10 & 20 & 10 & 7 & 20 & 6 & 4 & 20 & 12 \\
		 & $\epsilon = 0.10$ & 1 & 16 & 3 & 10 & 12 & 14 & 11 & 12 & 12 & 10 & 20 & 10 & 7 & 20 & 6 & 4 & 20 & 12 \\
		 & $\epsilon = \infty$ & 1 & 16 & 2 & 11 & 12 & 14 & 11 & 12 & 12 & 10 & 20 & 10 & 7 & 20 & 6 & 4 & 20 & 12 \\
		\hline & $\mathcal{K} = 1$ & 8 & 12 & 12 & 11 & 12 & 14 & 11 & 12 & 12 & 10 & 13 & 13 & 12 & 11 & 6 & 7 & 12 & 12 \\
		$*$ & $\mathcal{K} = 6$ & 1 & 16 & 4 & 8 & 12 & 14 & 11 & 12 & 12 & 10 & 20 & 10 & 7 & 20 & 6 & 5 & 20 & 12 \\
		 & $\mathcal{K} = 12$ & 19 & 11 & 8 & 2 & 12 & 14 & 10 & 1 & 2 & 1 & 20 & 6 & 6 & 20 & 20 & 5 & 20 & 12 \\
		 & $\mathcal{K} = 18$ & 1 & 11 & 7 & 8 & 12 & 14 & 11 & 1 & 13 & 4 & 20 & 8 & 16 & 20 & 10 & 5 & 20 & 12 \\
		 & $\mathcal{K} = 24$ & 13 & 11 & 8 & 1 & 12 & 14 & 1 & 11 & 3 & 1 & 20 & 17 & 20 & 20 & 6 & 4 & 20 & 11 \\
		 & $\mathcal{K} = 30$ & 12 & 11 & 14 & 7 & 12 & 14 & 11 & 5 & 11 & 9 & 8 & 20 & 7 & 5 & 8 & 20 & 6 & 12 \\
		\hline & $\beta = 0.5$ & 10 & 20 & 11 & 2 & 9 & 13 & 11 & 1 & 1 & 2 & 20 & 6 & 3 & 18 & 18 & 5 & 20 & 12 \\
		$*$ & $\beta = 1.0$ & 1 & 16 & 4 & 8 & 12 & 14 & 11 & 12 & 12 & 10 & 20 & 10 & 7 & 20 & 6 & 5 & 20 & 12 \\
		 & $\beta = 1.5$ & 1 & 17 & 5 & 10 & 12 & 14 & 11 & 12 & 12 & 10 & 16 & 10 & 7 & 20 & 6 & 5 & 20 & 12 \\
		 & $\beta = 2.0$ & 1 & 20 & 5 & 10 & 12 & 14 & 11 & 12 & 12 & 10 & 12 & 10 & 7 & 20 & 6 & 6 & 20 & 12 \\
		\hline
	\end{tabular}
	\caption{Solution vectors for all artificial trial sets as well as the initial solution vector. Rows listing a parameter $\epsilon$, $\mathcal{K}$, or $\beta$ correspond to the trial set in which the listed parameter was varied. The three starred ($*$) rows represent a single trial with $\epsilon = 0.01$, $\mathcal{K} = 6$, and $\beta = 1.0$.}
	\label{table:smallscalesol}
\end{table}

\end{document}